\documentclass[preprint,
3p,
times loads txfonts
]{elsarticle}

\usepackage{amsmath}

\usepackage{amsthm}
\usepackage{amssymb}
\usepackage{color}
\usepackage{xcolor}
\usepackage{cases}
\usepackage{cancel}
\usepackage{enumitem}% Improved itemize & enumerate environments.    
\usepackage{algorithm}
\usepackage{algpseudocode}
\usepackage{mathtools}
\usepackage[OMLmathsfit]{isomath}  % Loads also fixmath.sty (upper case italic greek letters, definition of \upDelta etc).
\usepackage{graphicx}
\usepackage{grffile}               % Fixes errors with images whose file names contain more than one dot, e.g., *.0000.png.
\usepackage[T1]{fontenc}
\usepackage{lmodern}
\usepackage[normalem]{ulem}
\usepackage{comment}
\usepackage{hyperref}

\usepackage[cmintegrals]{newpxmath}  % The fitting math symbols for newpxtext.
\usepackage{bm}                      % Provides \boldsymbol that we use in \vec{}.
\usepackage{subfiles}
\usepackage{float}
\usepackage{subfig}
\usepackage{epsfig}
\usepackage{psfrag}
\usepackage{graphics}
\usepackage[running]{lineno} % Show line numbers for review.
\newcommand*\patchAmsMathEnvironmentForLineno[1]{%
  \expandafter\let\csname old#1\expandafter\endcsname\csname #1\endcsname
  \expandafter\let\csname oldend#1\expandafter\endcsname\csname end#1\endcsname
  \renewenvironment{#1}%
     {\linenomath\csname old#1\endcsname}%
     {\csname oldend#1\endcsname\endlinenomath}}% 
\newcommand*\patchBothAmsMathEnvironmentsForLineno[1]{%
  \patchAmsMathEnvironmentForLineno{#1}%
  \patchAmsMathEnvironmentForLineno{#1*}}%
\AtBeginDocument{%
\patchBothAmsMathEnvironmentsForLineno{equation}%
\patchBothAmsMathEnvironmentsForLineno{align}%
\patchBothAmsMathEnvironmentsForLineno{flalign}%
\patchBothAmsMathEnvironmentsForLineno{alignat}%
\patchBothAmsMathEnvironmentsForLineno{gather}%
\patchBothAmsMathEnvironmentsForLineno{multline}%
}

\journal{{\color{red}{Not yet submitted.}}}

\begin{document}
% \linenumbers

\begin{frontmatter}

%%%%%%%%%%%%%%%%%%%%%%%%%%%%%%%%%%%%%%%%%%%%%%%%%%%%%%%
%
%   TITLE 
%
%%%%%%%%%%%%%%%%%%%%%%%%%%%%%%%%%%%%%%%%%%%%%%%%%%%%%%%
\title{A Structure-Preserving Penalization Method for the Single-species Rosenbluth-Fokker-Planck Equation}

\author[1,2]{Hamad El Kahza\corref{cor1}}\ead{helkahza@udel.edu}
\author[2]{Luis Chac{\'o}n}
\author[2]{William Taitano}
\author[1]{Jing-Mei Qiu}
\author[3]{Jingwei Hu}

\address[1]{Department of Mathematical Sciences, University of Delaware, Newark, DE, 19716, USA}
\address[2]{Theoretical Division, Los Alamos National Laboratory, Los Alamos, NM, 87545, USA}
\address[3]{Department of Applied Mathematics, University of Washington, Seattle, WA, 98195, USA}
\cortext[cor1]{Corresponding Author}
%%%%%%%%%%%%%%%%%%%%%%%%%%%%%%%%%%%%%%%%%%%%%%%%%%
%
%   Abstract
%
%%%%%%%%%%%%%%%%%%%%%%%%%%%%%%%%%%%%%%%%%%%%%%%%%%
%
%
%
%
\begin{abstract}
{The Rosenbluth-Fokker-Planck (RFP) equation describes Coulomb collisional dynamics within and across species in plasmas. It belongs to the broader class of anisotropic-diffusion-advection equations, whose numerical approximation is highly-nontrivial due to its nonlinearity, stiffness, and structural properties such as conservation and entropy dissipation (hence with the Maxwellian distribution as the equilibrium state). In this paper, we propose a structure-preserving penalization scheme for the stiff, single-species RFP equation. The scheme features three novel components: 1) a novel generalization of the well-known Chang–Cooper discretization for the RFP equation that is equilibrium-preserving and enables positivity while preserving mass, momentum, and energy; 2) an easy-to-invert isotropic variable-coefficient penalization operator to deal with the temporal stiffness without resorting to a fully implicit scheme, borrowing ideas from explicit-implicit-null (EIN) methods, and 3) an adaptive timestepping strategy that preserves the positivity of the full penalized scheme. 
The resulting scheme conserves mass, momentum, and energy strictly, is unconditionally stable, and robustly positivity preserving. The scheme is demonstrated with linear and nonlinear anisotropic diffusion examples of increasing complexity, including several single-species RFP examples.}
\end{abstract}

\begin{keyword}
Rosenbluth Fokker Planck \sep Explicit-Implicit-Null Penalization \sep Structure Preservation \sep Asymptotic Preservation \sep Positivity Preservation \sep Equilibrium Preservation \sep Anisotropic Diffusion
\vspace{.5\baselineskip}
%\MSC 65M12 \sep 65M60 
\end{keyword}
\end{frontmatter}
%%%%%%%%%%%%%%%%%%%%%%%%%%%%%%%%%%%%%%%%%%%%%%%%%%
%
%   Introduction
%
%%%%%%%%%%%%%%%%%%%%%%%%%%%%%%%%%%%%%%%%%%%%%%%%%%
\section{Introduction}
\label{sec:introduction}

We propose a novel structure-preserving penalization scheme for the single-species Rosenbluth-Fokker-Planck (RFP) equation that strictly conserves mass, momentum, and energy while preserving positivity and the analytical equilibrium exactly. The RFP equation is a nonlinear advection-diffusion-type model that describes Coulomb grazing-angle collisions in a plasma. It is mathematically equivalent to the well-known Landau-Fokker-Planck (LFP) equation \cite{landau1936kinetische}, but is more amenable to fast numerical integration, with computational complexities scaling linearly with the number of unknowns \cite{taitano2015mass} instead of quadratically.  

The RFP equation (as its analogous, the LFP) features rich mathematical structure. It strictly conserves mass, momentum, and energy, satisfies an H-theorem (strict growth of physical entropy except at equilibrium), and has the Maxwellian distribution in its kernel. Discretizing the RFP equation while preserving all this structure is non-trivial. Temporally, the RFP equation can be very stiff in moderately and strongly collisional regimes, benefiting from implicit timestepping techniques. However, these typically result in a coupled set of nonlinear algebraic equations of difficult inversion. Spatially, a suitable discretization should inherit strict conservation properties, the H-theorem, a maximum principle (to enable strict preservation of positivity), and strict equilibrium preservation (meaning that the analytical Maxwellian is in the kernel of the discretized RFP operator).

Structurally, the RFP is a nonlinear anisotropic-diffusion-advection operator. The diffusion anisotropy can prevent a discrete maximum principle regardless of timestepping strategy \cite{sharma2007preserving}, which is a necessary condition for positivity preservation. In equilibrium, advective and diffusive fluxes exactly balance for {\em any} Maxwellian distribution, implying that the RFP collision operator is infinitely degenerate. The degeneracy is resolved by the strict conservation of mass, momentum, and energy, which uniquely determine the steady state Maxwellian for a given initial condition.

Preserving all these structural properties discretely is problematic. Mass conservation is straightforward to preserve by discretizing the divergence-form of the operator conservatively. However, strict numerical conservation of momentum and energy in the RFP equation is difficult to enforce. This is unlike its cousin the LFP equation, for which conservation follows from operator symmetries easy to build in discretely. Early work on strict energy conservation for the RFP equation \cite{chacon2000discret,chacon2000implicit} was only partially successful owing to boundary effects. A complete solution for strict discrete conservation for the RFP was proposed in \cite{taitano2015mass}. The authors proposed a Lagrange-multiplier technique in which certain parameters are introduced in the discretization  of the RFP operator to guarantee strict conservation of mass, momentum, and energy. The approach works regardless of temporal discretization strategy, including fully implicit methods. This is the approach we will follow in this study.
 
Strict positivity and equilibrium preservation remain an open problem for the RFP equation. Na\"ive spatial discretizations of the RFP operator may cause severe positivity issues due to the lack of a maximum principle. The diffusion is of a variable anisotropic nature, meaning that, on a given computational mesh, the principal diffusion directions at each point in the domain need not align with the mesh. If this misalignment is not handled carefully, it can lead to a significant loss of positivity. Classical discretization approaches for anisotropic diffusion, such as for the anisotropic heat equation, rely on high-order finite differences for the anisotropic PDE components~\cite{deka2022exponential,crouseilles2015comparison}. 

Such discretizations typically lack a maximum principle, which may lead to nonphysical negativity in the solution unless extremely fine meshes are utilized. More involved, but less practical, strategies in the literature include non-negative directional splitting: the operator is locally rotated at the continuum level to align with the principal components of the diffusion tensor, and a wide stencil is then used to discretize along the rotated axes aligned with those principal components; see \cite{ngo2016monotone} and references therein. While this approach is attractive because the operator is rotated to align with the principal components of the diffusion tensor, under extreme anisotropy one must use a very wide differencing stencil to be able to discretize the rotated operator on the mesh. Moreover, the conservative form of the operator is broken, spoiling mass conservation.

More practical approaches, which we adopt in this work, ``advectionalize'' (a mathematical manipulation to be defined precisely later in this study) the diffusion fluxes that do not align with the mesh and treat them with a monotonicity-preserving spatial discretization\cite{du2018positivity}.
While advectionalization is attractive for its simplicity, it leads by construction to a strongly nonlinear discretization technique. Moreover, it is unclear {\em a priori} how to reconcile this discretization strategy with other requirements such as equilibrium preservation and strict conservation.
In this work, we resolve this tension by adapting (to our knowledge, for the first time) the equilibrium- and maximum-principle-preserving Chang-Cooper discretization to the RFP equation~\cite{chang1970practical}. 
Their technique (developed for a simplified Fokker-Planck equation, the Lenard-Bernstein collisional model \cite{lenard1958plasma}) consists of choosing the interpolation weights of the advective flux such that, at equilibrium, the numerical flux is exactly zero for the analytical Maxwellian. However, when adopted for the RFP setting, the equilibrium is not analytically preserved due to the numerical error introduced in the computation of the diffusion tensor and advection coefficients. We provide a modified Chang-Cooper scheme for the RFP equation that resolves this problem and leverages the advectionalized form of the RFP to produce a maximum-principle- and equilibrium-preserving discrete solution.  

Temporally, the discretized RFP equation is nonlinear and stiff. Explicit schemes are very limiting owing to numerical stability constraints in the timestep, which scale quadratically with spatial resolution. Implicit methods are desirable to circumvent such numerical stability constraints. However, by construction they lead to large systems of nonlinear algebraic equations, difficult to invert. Practical implicit schemes that conserve mass, momentum, and energy have been developed \cite{taitano2015mass}, but they require modern nonlinear iterative methods and scalable preconditioning (e.g., with multigrid methods), and in practice conservation properties are only preserved to the nonlinear tolerance employed. Moreover, the nonlocal nature of the RFP equation (via the Rosenbluth potentials) make usual preconditioning approximations brittle for sufficiently large timesteps unless sophisticated multiscale methods are used (such as high-order/low-order iterative schemes \cite{chacon2017multiscale}). 

To mitigate these temporal difficulties, so-called explicit-implicit-null (EIN) approaches have been developed with applications to diffusion~\cite{deka2022exponential,tan2025high,crouseilles2015comparison}, classical Boltzmann \cite{jin2011class}, and classical \cite{yan2013successive} and quantum \cite{hu2012numerical} Landau-Fokker-Planck equations. The EIN strategy treats the original (typically hard to invert) nonlinear operator explicitly, while adding and subtracting a carefully crafted linear operator that is handled in an IMEX fashion. The ``null'' in EIN expresses the need of the penalization operator to respect the null space of the original system, when present. When properly chosen, the penalization operator unconditionally stabilizes the scheme, allowing arbitrary timesteps (albeit constrained by accuracy). At each time step, only one linear penalization operator needs to be inverted, which is computationally attractive. Directly related to our study is the work in \cite{yan2013successive}, where the linear Lenard-Bernstein-Fokker-Planck (LBFP) operator is used to penalize the nonlinear RFP equation. LBFP is null-space preserving, and shown to be much more accurate than simpler penalization operators such as the Laplacian. However, LBFP still features constant coefficients, which have been recently shown~\cite{tan2025high} to degrade accuracy due to overpenalization in key areas of the domain. In the reference, a variable-coefficient penalization approach is proposed and shown to be much more accurate. We adopt this approach here.

For the penalization of the RFP equation, to ensure accuracy, it is therefore instrumental to devise a null-space preserving penalization operator with a variable penalization coefficient. Here, we explore a minor reformulation of the LBFP operator that introduces a variable isotropic penalization coefficient while still being null-space preserving. The proposed penalization operator can be made strictly conservative of mass, momentum, and energy by introducing a mild nonlinear iteration for effective drift and temperature parameters that converge to the exact Maxwellian ones in equilibrium. A Chang-Cooper discretization of the so-formulated penalization operators ensures the existence of a discrete maximum principle, a necessary condition for positivity preservation, and preservation of the analytical null space.
Strict positivity preservation for the full discrete penalized scheme is achieved here with a final ingredient: a careful timestep choreography. Following \cite{barnett2004fast}, an adaptive timestep strategy, derived from a Fourier analysis of the penalized scheme, is proposed that features exceptionally robust positivity preservation properties.

In summary, this work combines anisotropy-aware spatial discretizations with a penalty-based temporal integrator for the RFP equation. We (i) construct a modified Chang-Cooper-type spatial discretization tailored to the RFP operator that analytically preserves the Maxwellian equilibrium when near steady state, while employing an “advectionalized’’ treatment of misaligned anisotropic diffusion; (ii) enforce exact conservation of mass, momentum, and energy through a Lagrange-multiplier modification of both the RFP and penalization operators with variable coefficients; and (iii) prescribe an adaptive time stepping, that enables large time steps while recovering positivity in practice. These ingredients form the basis of the conservative, positivity preserving, and steady-state preserving penalized solver for the RFP equation proposed and analyzed in this paper.

The rest of the paper is organized as follows. In Section~\ref{sec:RFP-equation}, we recall the Rosenbluth-Fokker-Planck equation and its conservation properties. Section~\ref{sec:EIN_anisotropic_heat} introduces the EIN penalization method for anisotropic diffusion, together with stability and positivity analyses and an adaptive positivity-preserving time-stepping strategy and spatial discretization. In Section~\ref{sec:EIN-RFP}, we develop a structure-preserving EIN penalization scheme for the single-species RFP operator, emphasizing conservative penalization and null-space preservation. Section~\ref{sec:RFP_spatial_discretization} presents a compatible spatial discretization of the EIN-penalized RFP equation that enforces a discrete maximum principle, preserves the null space, and yields a conservative nonlinear solution strategy. Numerical results in Section~\ref{sec:numerical_results} verify the method with anisotropic diffusion, cross-species pitch-angle scattering, and nonlinear RFP relaxation tests. Finally, Section~\ref{sec:conclusions} offers concluding remarks and directions for future work.

%%%%%%%%%%%%%%%%%%%%%%%%%%%%%%%%%%%%%%%%%%%%%%%%%%
%
%   Rosenbluth-Fokker-Planck equation
%
%%%%%%%%%%%%%%%%%%%%%%%%%%%%%%%%%%%%%%%%%%%%%%%%%%
\section{The single-species Rosenbluth-Fokker-Planck equation}
\label{sec:RFP-equation}
The Rosenbluth-Fokker-Planck (RFP) collision operator is an exact reformulation of the Landau-Fokker-Planck operator, and can be written as \cite{rosenbluth1957fokker}:
\begin{equation}
\label{eq:RFP}
    \partial_t f
    \;=\;
    \mathcal{C}(f)
    :=
  \nabla_\mathbf{v} \cdot\bigl( {\cal D} \cdot \nabla_\mathbf{v} f - \mathbf{A} f\bigr),
\end{equation}
where $\nabla_\mathbf{v}$ denotes the gradient in velocity space.
The diffusion tensor ${\cal D}$ and friction vector $\mathbf{A}$ are
defined in terms of the Rosenbluth potentials $G$ and $H$ as
\begin{equation}
\label{eq:DA_def}
\left\{
\begin{aligned}
{\cal D} &= \nabla_\mathbf{v} \nabla_\mathbf{v} G,\\
\mathbf{A} &= \nabla_\mathbf{v} H,
\end{aligned}
\right.
\end{equation}
where the potentials $G$ and $H$ satisfy the coupled Poisson equations
\begin{equation}
    \label{eq:rosenbluth_poisson}
    \nabla_\mathbf{v}^2 H = -8\pi f \,\,,\,\,
    \nabla_\mathbf{v}^2 G = H.
\end{equation}
Here, $\nabla_\mathbf{v} \nabla_\mathbf{v}$ is the Hessian, and $\nabla_\mathbf{v}^2=\nabla_\mathbf{v} \cdot \nabla_\mathbf{v}$ is the Laplacian operator.

The RFP operator rigorously conserves mass, momentum, and energy, which can be
expressed as:
\begin{flalign}
    \label{eq:conservation_symmetry}
    &
    \left<
        \boldsymbol{\phi},
        \nabla_\mathbf{v} \cdot
        \left[
         {\cal D} \cdot \nabla_\mathbf{v} f
        -
        \mathbf{A}f
        \right]
    \right>_v = \mathbf{0},
\end{flalign}
where 
\[
    \boldsymbol{\phi} = \bigl(1,\, \mathbf{v},\, |\mathbf{v}|^2\bigr)^\top
\]
and 
\[
    \left< f(\mathbf{v}),g(\mathbf{v}) \right>_v 
    = \int_{\mathbb{R}^3} d\mathbf{v}\, f(\mathbf{v})g(\mathbf{v})
\]
denotes the velocity-space inner product. 
{Furthermore, the RFP operator features an H-theorem \cite{hinton1983collisional}, which dictates that entropy increases monotonically except at equilibrium. The equilibrium solutions are those in the kernel of the RFP operator, which are Maxwellians. Owing to the conservation properties stated above, the RFP operator admits a unique Maxwellian kernel for a given initial condition, parametrized by the conserved density $n_M$, drift velocity $\mathbf{u}^M$, and thermal velocity $v_{th,M}$. In $d$-dimensional velocity space (i.e., $\mathbf{v}\in\mathbb{R}^d$), this Maxwellian is given by
\begin{equation}
\label{Maxwellian_kernel}
    f^M(\mathbf{v})
    =
    \frac{n_M}{(2\pi)^{d/2}\,v_{th,M}^d}
    \exp\!\left(
        -\frac{(\mathbf{v}-\mathbf{u}^M)^2}{2\,v_{th,M}^2}
    \right).
\end{equation}
Here, the parameters are defined through the mass of the colliding species, $m$, and conserved moments as
\begin{equation}
    v_{th,M} = \sqrt{\frac{T}{m}}, 
    \qquad
    n_M = \left< 1, f^M \right>_v,
    \qquad
    n_M \mathbf{u}^M = \left< \mathbf{v}, f^M \right>_v,
    \qquad
    d\,n_M T = m \left< |\mathbf{v}-\mathbf{u}^M|^2, f^M \right>_v .
\end{equation}
}

In this study, without loss of generality for the development of a penalization operator, we
restrict attention to distributions that are azimuthally symmetric in
velocity space, so that all quantities depend only on the cylindrical
coordinates $(v_\perp,v_\parallel)$ and are invariant with respect to the
angle~$\phi$. Accordingly, $d\mathbf{v}=2\pi v_\perp dv_\perp dv_\parallel$ and the inner product
reduces to
\[
    \left< f,g \right>_v
    =
    2\pi
    \int_{-\infty}^{\infty} dv_\parallel
    \int_{0}^{\infty} dv_\perp\, v_\perp\,
    f(v_\perp,v_\parallel)\,g(v_\perp,v_\parallel).
\]
The velocity-space Laplacian acting on a scalar $f(v_\perp,v_\parallel)$ is given by:
\[
\nabla^2_{\mathbf{v}} f
=
\frac{1}{v_\perp}\frac{\partial}{\partial v_\perp}\!\left(
v_\perp\,\frac{\partial f}{\partial v_\perp}
\right)
+
\frac{\partial^2 f}{\partial v_\parallel^2},
\]
and the Hessian by:
\[
\nabla_{\mathbf{v}}\nabla_{\mathbf{v}} f
=
\begin{pmatrix}
\displaystyle \frac{\partial^2 f}{\partial v_\perp^2}
&
\displaystyle \frac{\partial^2 f}{\partial v_\perp\,\partial v_\parallel}
\\[1.2ex]
\displaystyle \frac{\partial^2 f}{\partial v_\parallel\,\partial v_\perp}
&
\displaystyle \frac{\partial^2 f}{\partial v_\parallel^2}
\end{pmatrix}.
\]

Our goal is to solve \eqref{eq:RFP} with a suitable EIN penalization method such that it is unconditionally stable, and preserves positivity, strict conservation of mass, momentum, and energy, and the analytical null space (Maxwellian distribution). 
In what follows, we consider the penalization operator for a linear anisotropic diffusion equation as a proxy for the nonlinear RFP advection-diffusion equation of interest. We discuss stability and positivity in this context first, to be generalized later to the penalized RFP equation.

%%%%%%%%%%%%%%%%%%%%%%%%%%%%%%%%%%%%%%%%%%%%%%%%%%
%
%   EIN Method for Anisotropic Diffusion
%
%%%%%%%%%%%%%%%%%%%%%%%%%%%%%%%%%%%%%%%%%%%%%%%%%%
\section{EIN penalization method for anisotropic diffusion}
\label{sec:EIN_anisotropic_heat}

We consider the generic anisotropic diffusion equation:
\begin{equation}
  \partial_t f \;=\; \nabla \cdot \bigl( \bar{\bar{D}}\cdot \,\nabla f \bigr)
  \;=:\; \mathcal{Q}(f),
  \label{eq:anisotropic_diffusion}
\end{equation}
posed on a rectangular spatial domain $\Omega = (a,b)\times(c,d)$ with suitable initial
and boundary conditions. The diffusion tensor $\bar{\bar{D}}(\mathbf{x})$ (with $\mathbf{x}=(x,y)$) is assumed symmetric
and uniformly positive definite on $\Omega$, and is given by:
\begin{equation}
\label{eq:diff_coef_def}
    \bar{\bar{D}}(\mathbf{x})
    \;=\;
    \begin{pmatrix}
      a(\mathbf{x}) & b(\mathbf{x}) \\
      b(\mathbf{x}) & c(\mathbf{x})
    \end{pmatrix}.
\end{equation}

We introduce a penalization operator in the spirit of the EIN approach
(e.g.\ \cite{jin2011class, yan2013successive, tan2025high}), chosen so that its discretization is linear and easy
to invert. Using a backward Euler discretization in time, the penalized
update takes the form:
\begin{equation}
    f^{n+1} - \Delta t\, \mathcal{L}_\beta(f^{n+1})
    \;=\;
    f^{n}
    + \Delta t\Bigl(\mathcal{Q}(f^{n}) - \mathcal{L}_\beta(f^{n})\Bigr),
\end{equation}
where $\mathcal{L}_\beta$ is a linear isotropic-diffusion operator that captures the stiffness
of the original operator $\mathcal{Q}$, thereby stabilizing the scheme,
but is easier to invert. In the spirit of~\cite{tan2025high}, we consider the variable-coefficient isotropic diffusion operator:
\begin{equation}
    \mathcal{L}_\beta(f)
    \;=\;
    \nabla \cdot \bigl( \beta(\mathbf{x})\,\nabla f \bigr),
\end{equation}
with a scalar coefficient $\beta(\mathbf{x})>0$ to be determined to ensure global stability while avoiding overpenalization. We discuss this next in a semi-discrete context (temporally discrete, spatially continuous), and later generalize it to the fully discrete scheme.

\subsection{Stability analysis of the semi-discrete scheme}

For the stability analysis, we consider a constant anisotropic diffusion tensor, which we eigen-decompose as:
\begin{equation}
    \bar{\bar{D}} \;=\;
    \begin{pmatrix}
        a & b \\
        b & c
    \end{pmatrix}
    \;=\;
    \lambda_1 \mathbf{v}_1 \mathbf{v}_1^\top \,+\, \lambda_2 \mathbf{v}_2 \mathbf{v}_2^\top,
\end{equation}
for some orthonormal basis $\{\mathbf{v}_1, \mathbf{v}_2\}$ and eigenvalues
$0<\lambda_2 \le \lambda_1$. Transforming the anisotropic diffusion equation into the rotated Cartesian coordinates $(\zeta,\eta)$ defined by
$(\mathbf{v}_1,\mathbf{v}_2)$, we obtain the transformed penalized update
\begin{equation}
    f^{n+1}
    - \Delta t\,\beta\bigl(\partial_{\zeta\zeta}^2 f^{n+1}
    + \partial_{\eta\eta}^2 f^{n+1}\bigr)
    \;=\;
    f^n
    + \Delta t\Bigl[(\lambda_1 - \beta)\,\partial_{\zeta\zeta}^2 f^n
    + (\lambda_2 - \beta)\,\partial_{\eta\eta}^2 f^n\Bigr].
\end{equation}
Applying a Fourier transform in $(\zeta,\eta)$, each Fourier mode with
wave numbers $\mathbf{k}=(k_1,k_2)$ satisfies
\begin{equation}
    \hat f^{\,n+1}(\mathbf{k})
    \;=\;
    R(\mathbf{k})\,\hat f^{\,n}(\mathbf{k}),
\end{equation}
where the amplification factor $R$ is given by
\begin{equation}
\label{eq:amp_fac_fourier}
    R(\mathbf{k})
    \;=\;
    \frac{1 - \Delta t\bigl((\lambda_1-\beta) k_1^2 + (\lambda_2-\beta) k_2^2\bigr)}
         {1 + \Delta t\,\beta\bigl(k_1^2 + k_2^2\bigr)}
    \;=\;
    1 - \Delta t\,
    \frac{\lambda_1 k_1^2 + \lambda_2 k_2^2}{1 + \Delta t\,\beta\bigl(k_1^2 + k_2^2\bigr)}.
\end{equation}
The stability condition is $|R(\mathbf{k})| \le 1$ for all
$(\mathbf{k})$. A
sufficient condition for stability, in terms of the largest eigenvalue $\lambda_1$ of $D$, is: 
\begin{equation}
    \beta \;\geq\; \frac{\lambda_1}{2}.
\end{equation}

For non-constant diffusion tensors, one option is to choose $\beta$ greater than half of the
largest eigenvalue of the diffusion tensor over the entire domain. 
This provides stability, but at the cost of accuracy due to over-penalization (see~\cite{tan2025high} and also the numerical results section). Instead, inspired by the reference, we introduce a variable-coefficient
penalization such that:
\begin{equation}
\label{eq:beta_def}
    \beta(\mathbf{x}) \;\geq\; \frac{\lambda_1(\mathbf{x})}{2}.
\end{equation}
This choice provides the right amount of penalization locally for stability while maintaining reasonable accuracy, as we will demonstrate numerically later in this study.

\subsection{Positivity analysis of the semi-discrete scheme}

The penalized temporal discretization of the anisotropic diffusion operator
is in general not positivity preserving. To see this, we express the
solution update as a convolution with the inverse Fourier transform of the
amplification factor $R(\mathbf{k})$:
\begin{equation}
    f^{n+1}
    \;=\;
    \mathcal{F}^{-1}\{R(\mathbf{k})\} * f^{n}.
    \label{solution_update_fourier}
\end{equation}
Assuming that the previous solution $f^n$ is strictly positive, the new
solution $f^{n+1}$ remains positive if the convolution kernel
\[
  R(\zeta,\eta) := \mathcal{F}^{-1}\{R(\mathbf{k})\}(\zeta,\eta)
\]
is nonnegative. 
Up to a normalization constant, the inverse Fourier
transform can be written as:
\begin{align}
    R(\zeta,\eta)
    &= \delta(\zeta,\eta)
    + \frac{1}{\beta}\Bigl[
        \lambda_1\,\partial_{\zeta\zeta}^2
        K_0\!\Bigl(\frac{\sqrt{\zeta^2+\eta^2}}{\sqrt{\Delta t\,\beta}}\Bigr)
        +
        \lambda_2\,\partial_{\eta\eta}^2
        K_0\!\Bigl(\frac{\sqrt{\zeta^2+\eta^2}}{\sqrt{\Delta t\,\beta}}\Bigr)
      \Bigr],
\end{align}
where $K_0$ is the zeroth-order modified Bessel function of the second kind and $\delta$ is the Dirac delta.
Separating the eigenvalues into
their mean and deviation,
\[
    \bar\lambda := \frac{\lambda_1 + \lambda_2}{2},
    \qquad
    \delta\lambda := \frac{\lambda_1 - \lambda_2}{2},
\]
we can write:
\[
    \lambda_1\,\partial_{\zeta\zeta}^2 + \lambda_2\,\partial_{\eta\eta}^2
    \;=\;
    \bar\lambda\,(\partial_{\zeta\zeta}^2 + \partial_{\eta\eta}^2)
    + \delta\lambda\,(\partial_{\zeta\zeta}^2 - \partial_{\eta\eta}^2).
\]
Noting the identity $\Delta K_0 = K_0$ for the Laplacian operator $\Delta$, we obtain (again up to a normalization factor):
\begin{align}
    R(\zeta,\eta)
    &= \delta(\zeta,\eta)
    + \frac{1}{\Delta t\,\beta^2}
    \Biggl[
        \bar\lambda\,
        K_0\!\Bigl(\frac{r}{\sqrt{\Delta t\,\beta}}\Bigr)
        +
        \delta\lambda\,
        \frac{\zeta^2 - \eta^2}{\zeta^2 + \eta^2}\,
        K_2\!\Bigl(\frac{r}{\sqrt{\Delta t\,\beta}}\Bigr)
    \Biggr],
\end{align}
where $r = \sqrt{\zeta^2 + \eta^2}$ and $K_2$ is the modified Bessel
function of order two.
Introducing polar coordinates:
\[
    \zeta = r\cos\theta, \qquad \eta = r\sin\theta,
\]
we obtain:
\begin{equation}
\label{eq:amp_factor_r_th}
    R(r,\theta)
    \;=\;
    \frac{\delta(r)}{r}
    + \frac{1}{\Delta t\,\beta^2}
    \Biggl[
      \bar\lambda\,
      K_0\!\left(\dfrac{r}{\sqrt{\beta\,\Delta t}}\right)
      +
      \underbrace{\delta\lambda\,\cos(2\theta)\,
      K_2\!\left(\dfrac{r}{\sqrt{\beta\,\Delta t}}\right)}_{\text{sign-changing}}
    \Biggr].
\end{equation}
The explicit angular dependence through $\cos(2\theta)$ is sign-changing,
so the kernel $R$ may become negative, and the scheme is
therefore not positivity preserving for arbitrarily large time steps.
However, the subsequent analysis shows that a time-adaptive strategy can be devised such that positivity is largely preserved. We discuss this next.

\subsection{Adaptive positivity-preserving time stepping}\label{Sec:Time_Adaptivity}

During the time update~\eqref{solution_update_fourier} in its convolution
representation, it is the anisotropic contribution in \eqref{eq:amp_factor_r_th},
\[
     \frac{1}{\Delta t\,\beta^2} \,\delta\lambda\,\cos(2\theta)\,
  K_2\!\left(\frac{r}{\sqrt{\beta\,\Delta t}}\right)
\]
that may induce negativity in $f^{n+1}$. In the isotropic case
$\lambda_1 = \lambda_2$, this term vanishes and the resulting scheme is
positivity preserving regardless of timestep. In the strongly anisotropic case, e.g., $\lambda_2 \ll \lambda_1$, the time step must be chosen carefully to
ensure that this anisotropic contribution does not introduce negativity.

In fact, an adaptive timestep strategy is possible that ensures positivity while allowing large timesteps. It is based on limiting the contribution of the sign-indefinite term in the amplification
factor, and it is best derived in its Fourier space expression, Eq. \eqref{eq:amp_fac_fourier}, which can be rewritten as:
\begin{equation}
\label{eq:amp_fac_fourier_2}
    R(\mathbf{k})
    \;=\;
    1 - \Delta t \bar{\lambda}\,
    \frac{k_1^2 + k_2^2}{1 + \Delta t\,\beta\bigl(k_1^2 + k_2^2\bigr)}
    - \Delta t\,\delta\lambda\,
    \frac{k_1^2 - k_2^2}{1 + \beta\,\Delta t\,(k_1^2 + k_2^2)}.
\end{equation}
The Fourier component corresponding to the sign-changing (anisotropic) part of the
amplification factor is:
\begin{equation}
    R_{\mathrm{aniso}}(\mathbf{k})
    \;=\;
    \Delta t\,\delta\lambda\,
    \frac{k_1^2 - k_2^2}{1 + \beta\,\Delta t\,(k_1^2 + k_2^2)}.
\end{equation}
At any given time $t$, we require this anisotropic contribution to remain
small in the sense that
\begin{equation}
    \bigl|R_{\mathrm{aniso}}(\mathbf{k})\bigr|
    \;<\; \varepsilon,
\end{equation}
for all relevant wave numbers $(\mathbf{k})$, where $\varepsilon>0$ is a
user-specified tolerance. This condition will drive our adaptive choice
of the time step $\Delta t$.

Introducing polar coordinates in Fourier space as:
\[
    |\mathbf{k}|^2 = k_1^2 + k_2^2, \qquad
    k_1 = |\mathbf{k}|\cos\sigma,\quad
    k_2 = |\mathbf{k}|\sin\sigma,
\]
we can write:
\[
    k_1^2 - k_2^2 = |\mathbf{k}|^2(\cos^2\sigma - \sin^2\sigma)
    = |\mathbf{k}|^2\cos(2\sigma),
\]
and hence
\begin{equation}
    R_{\mathrm{aniso}}(\mathbf{k})
    = \Delta t\,\delta\lambda\,
    \frac{|\mathbf{k}|^2\cos(2\sigma)}{1 + \beta\,\Delta t\,|\mathbf{k}|^2}.
\end{equation}
Taking absolute values and using $|\cos(2\sigma)|\le 1$ gives
\begin{equation}
    \bigl|R_{\mathrm{aniso}}(\mathbf{k})\bigr|
    \;\le\;
    \Delta t\,|\delta\lambda|\,
    \frac{|\mathbf{k}|^2}{1 + \beta\,\Delta t\,|\mathbf{k}|^2}
    \;\le\;
    \Delta t\,|\delta\lambda|\,|\mathbf{k}|^2.
\end{equation}

As the solution evolves under diffusion, high-frequency modes are rapidly
damped, and at time $t$ only a portion of the spectrum such that only modes satisfying:
\[
    |\mathbf{k}|^2 = O\!\left(\frac{1}{\beta t}\right)
\]
remain dynamically relevant~\cite{barnett2004fast}. 
Restricting our timestep to resolve only undamped modes, we obtain the
following timestep limiting condition:
\[
    \bigl|R_{\mathrm{aniso}}(\mathbf{k})\bigr|
    \;\lesssim\;
    \Delta t\,\frac{|\delta\lambda|}{\beta}\,\frac{1}{t}
    \;<\; \varepsilon,
\]
which in turn yields the time-step restriction:
\begin{equation}
    \Delta t
    \;\lesssim\;
    \varepsilon\,\frac{\beta}{|\delta\lambda|}\,t
    =
    \frac{2\,\varepsilon\,\beta}{\lambda_1 - \lambda_2}\,t,
    \label{eq:dt-bound-aniso}
\end{equation}
where we used $\delta\lambda = (\lambda_1-\lambda_2)/2$ and
$\lambda_1\ge\lambda_2$ by ansatz.
Allowing for the possibility of temporal change in the tensor-diffusion coefficient (as will be the case for the RFP equation), and letting $\alpha_n=\frac{2\,\varepsilon\,\beta_n}{\lambda_{1,n}-\lambda_{2,n}}$, where the subscript $n$ indicates temporal level, we obtain
a logarithmic time-stepping of the form
$\Delta t_n = \alpha_n t_n$. Solving for the time step,
we obtain:
\begin{align}
    \Delta t_n
    &= t_n - t_{n-1}= \frac{1}{\alpha_n}\,\Delta t_n
       - \frac{1}{\alpha_{n-1}}\,\Delta t_{n-1},\\
    \Rightarrow\quad
    \Delta t_n
    &= \frac{\alpha_n}{1-\alpha_n}\,
       \frac{1}{\alpha_{n-1}}\,\Delta t_{n-1},
\end{align}
for $n\geq 2$. In the fully discrete setting, we initialize the timestep with the explicit stability constraint for a given spatial discretization. Note that $\alpha_n<1$ for well-posedness of this formula, yielding an upper bound on the
threshold $\varepsilon$:
\[
    \varepsilon < \frac{\lambda_1-\lambda_2}{2\,\beta}.
\]
As $\varepsilon \rightarrow 0$, the timestep growth is negligible and we recover explicit time-stepping. When $\varepsilon \rightarrow \frac{\lambda_1-\lambda_2}{2\beta}$, then $\alpha \rightarrow 1$, and the prescription allows arbitrarily large timesteps, but this may result in solution negativity. One must therefore judiciously choose $\varepsilon$ to sufficiently damp the modes responsible for negativity in the anisotropic diffusion tensor while still achieving enough exponential growth in the time steps to benefit from sufficiently large timesteps. In our simulations, we set
\[
    \varepsilon = \min\!\left(0.05, \frac{\lambda_1-\lambda_2}{4\beta}\right),
\]
which we have found effective for all simulations in this study. 
We note that the upper bound \(0.05\) may be problem dependent and is selected here based on numerical testing. Tighter values typically lead to a slower growth of the time step, whereas larger values may induce negativity. 

\subsection{Maximum-principle-preserving spatial discretization}\label{MP_discretization_Heat}

The analysis so far has used a semi-discrete scheme, with continuous space and discrete time. While stability is expected to survive spatial discretization for any reasonable choice, the maximum principle will not. In fact, the anisotropic heat transport equation is well known to lack a maximum principle even when explicit stability timestep constraints are respected \cite{sharma2007preserving}.

In what follows, we follow \cite{du2018positivity} for a spatial discretization of the anisotropic transport equation with a maximum principle. The key idea is to discretize diagonal terms in the tensor diffusion equation with the usual three-point second-order approximations, which are known to feature a maximum principle, and off-diagonal terms as an advective operator after suitable reformulation (which we term here ``advectionalization''). The resulting discretization is nonlinear, but features a maximum principle.

We begin by rewriting the anisotropic diffusion operator in flux form,
\begin{equation}
   \mathcal{Q}(f)
    \;=\;
    \nabla\!\cdot\bigl(\bar{\bar{D}} \cdot \nabla f\bigr)
    \;=\;
    \partial_{x} J_x
    + \partial_{y} J_y,
\end{equation}
with flux $\mathbf{J}=(J_x,J_y)^\top = \bar{\bar{D}} \cdot \nabla f$, with $\bar{\bar{D}}$ defined in \eqref{eq:diff_coef_def}. 
Specializing for two dimensions without loss of generality, we split the flux in the $x$- and $y$-directions into their diagonal
(``isotropic'') and off-diagonal (``anisotropic'') contributions,
$\mathbf{J}=\mathbf{J}^{\text{iso}}+\mathbf{J}^{\text{aniso}}$, with
\begin{align}
    J_x^{\mathrm{iso}} &= a(\mathbf{x})\,\partial_x f,
    &
    J_y^{\mathrm{iso}} &= c(\mathbf{x})\,\partial_y f,
    \\
    J_x^{\mathrm{aniso}} &= b(\mathbf{x})\,\partial_y f,
    &
    J_y^{\mathrm{aniso}} &= b(\mathbf{x})\,\partial_x f.
\end{align}
To ``advectionalize'' the anisotropic part, we rewrite the off-diagonal
fluxes in an advective form by multiplying and dividing by $f$ (which is allowed because $f>0$ everywhere):
\begin{align}
    J_x^{\mathrm{aniso}}
    &= b(\mathbf{x})\,\partial_y f
     = b(\mathbf{x})\,\bigl(\partial_y \ln f\bigr)\,f,
    \\
    J_y^{\mathrm{aniso}}
    &= b(\mathbf{x})\,\partial_x f
     = b(\mathbf{x})\,\bigl(\partial_x \ln f\bigr)\,f.
\end{align}
This yields the advection-diffusion operator:
\begin{equation}
\label{eq:adv_diff_op}
    \mathcal{Q}(f)
    =
    \partial_x\bigl(a(\mathbf{x})\,\partial_x f\bigr)
    +
    \partial_y\bigl(c(\mathbf{x})\,\partial_y f\bigr)
    +
    \partial_x\bigl(\phi_x\,f\bigr)
    +
    \partial_y\bigl(\phi_y\,f\bigr),
\end{equation}
where the advection coefficients are given by:
\[
    \phi_x(\mathbf{x}) = b(\mathbf{x})\,\partial_y \ln f,
    \qquad
    \phi_y(\mathbf{x}) = b(\mathbf{x})\,\partial_x \ln f.
\]
Equation \eqref{eq:adv_diff_op} is discretized as follows. Diagonal (diffusive) fluxes are discretized
using a standard second-order finite-difference approximation, while the advective
operators are discretized with a suitable high-order monotonicity-preserving scheme. In this study, for the generic anisotropic diffusion tests, we use a variant of
the SMART scheme~\cite{gaskell1988curvature} for its simplicity and robustness in strongly anisotropic contexts \cite{chacon2025robust}.

%%%%%%%%%%%%%%%%%%%%%%%%%%%%%%%%%%%%%%%%%%%%%%%%%%
%
%   Application to Rosenbluth-Fokker-Planck
%
%%%%%%%%%%%%%%%%%%%%%%%%%%%%%%%%%%%%%%%%%%%%%%%%%%
\section{Structure-preserving EIN penalization scheme for the single-species RFP operator}
\label{sec:EIN-RFP}

The single-species RFP operator is substantially more complicated than the simple anisotropic-diffusion equation that we have analyzed so far. It is nonlinear, strictly conservative of mass, momentum, and energy, possesses a maximum principle, and features an infinitely degenerate null space (any Maxwellian is in its kernel).
A suitably penalized temporal and velocity-space discretization strategy should provide discrete solutions to maintain all these structural properties. However, to our knowledge, no such solution exists.
In the following sections, we build on the earlier developments to provide a comprehensive structure-preserving discrete penalization strategy for the RFP equation. The key ingredients are:
\begin{itemize}
    \item A structure-preserving semi-discrete variable-coefficient EIN penalization operator for the RFP equation based on a variation of the Lenard-Bernstein (LB) operator. The variable coefficient feature will result in significant accuracy gains vs. the constant-coefficient version proposed in \cite{jin2011class}, as will be demonstrated in the results section. The proposed penalization semi-discrete scheme is unconditionally stable, strictly conservative and null-space-preserving.
    \item A suitable spatial discretization for both the penalization and RFP operators that retains strict conservation of mass, momentum, and energy, preservation of the analytical null space, and a maximum principle.
\end{itemize}
We describe in this section the variable EIN penalization in the semi-discrete context (temporally discrete, spatially continuous). We will introduce the spatial discretization details later in this study.

\subsection{Conservative penalization approach}
\label{sec:cons_penalization}

We penalize the RFP operator with a variable-coefficient advection-diffusion equation (with isotropic diffusion) given by:
\begin{equation}
\label{eq:rfp_penalization_op}
    \mathcal{L}_\beta(f)
    =
    \nabla_\mathbf{v}\!\cdot\Bigl(\beta(\mathbf{v})\,
       \bigl[\nabla_\mathbf{v} f + \tfrac{\mathbf{v}-\mathbf{u}_\beta}{v_{th,\beta}^2}f\bigr]\Bigr),
\end{equation}
where, for azimuthally symmetric cylindrical coordinate systems, $\mathbf{u}_\beta=(0,u_{\parallel,\beta})^\top$, $ \mathbf{v} = (v_\perp, v_\parallel)^\top$, and $v_{th,\beta}=\sqrt{T/m}$. 
The time-discrete penalized RFP equation then reads:
\begin{equation}
\label{eq:penalized_RFP}
    f^{n+1}_{i,j}
    - \Delta t\,(\mathcal{L}_\beta f^{n+1})_{i,j}
    =
    f^{n}_{i,j}
    - \Delta t\,(\mathcal{L}_\beta f^{n})_{i,j}
    + \Delta t\,\mathcal{C}(f^{n})_{i,j}.
\end{equation}
The flux in \eqref{eq:rfp_penalization_op} is inspired in the well-known Lenard-Bernstein collision operator, but weighted with a variable $\beta(v_\perp , v_\parallel)$ isotropic diffusion coefficient (determined from stability as in \eqref{eq:beta_def}). 

The parameters
$u_{\parallel,\beta}$ and $v_{th,\beta}$ are moments of the $\beta$-weighted distribution function $f$, defined so as to enforce strict conservation of momentum and energy. This is done as follows.
We consider the moments vector 
\[
    \boldsymbol{\phi}(\mathbf{v}) = \bigl(1,\, \mathbf{v},\, v^2\bigr)^\top
\]
with $v^2=\mathbf{v}\cdot\mathbf{v}$,  and integrate the penalized time update against these moments:
\begin{equation}
    \int_\Omega
    \boldsymbol{\phi}(\mathbf{v})\,
    \Bigl(
        f^{n+1}
        - \Delta t\,(\mathcal{L}_\beta f^{n+1})
    \Bigr)\,d\mathbf{v}
    =
    \int_\Omega
    \boldsymbol{\phi}(\mathbf{v})\,
    \Bigl(
        f^{n}
        - \Delta t\,(\mathcal{L}_\beta f^{n})
        + \Delta t\,\mathcal{C}(f^{n})
    \Bigr)\,d\mathbf{v}.
\end{equation}
Since the discrete Rosenbluth-Fokker-Planck
collision operator $\mathcal{C}$ satisfies:
\begin{equation}
    \int_\Omega
    \boldsymbol{\phi}(\mathbf{v})\,
    \mathcal{C}(f^{n})\,d\mathbf{v} = \mathbf{0},
\end{equation}
and since
mass conservation is already enforced by the divergence form of the penalization
operators, strict momentum and energy conservation in \eqref{eq:penalized_RFP} simply requires:
\begin{equation}
\label{eq:lb_cons_constraints}
    \int_\Omega
    \begin{pmatrix}
        \mathbf{v} \\[0.2em]
        v^2
    \end{pmatrix}
    (\mathcal{L}_\beta f^{n+1})\,d\mathbf{v}
    = \mathbf{0},
    \qquad
    \int_\Omega
    \begin{pmatrix}
        \mathbf{v} \\[0.2em]
        v^2
    \end{pmatrix}
    (\mathcal{L}_\beta f^{n})\,d\mathbf{v}
    = \mathbf{0}.
\end{equation}
The parameters $(\mathbf{u}_\beta^{n+1},\lambda_\beta^{n+1})$ (with $\lambda = v_{th}^2$) and
$(\mathbf{u}_\beta^{n},\lambda_\beta^{n})$ are chosen so that
$\mathcal{L}_\beta^{n+1}$ and $\mathcal{L}_\beta^{n}$ are each conservative
in the second and third moments according to \eqref{eq:lb_cons_constraints}.
The conservation constraints lead (after integration by parts) to the following coupled equations:
\begin{equation}
    \lambda_\beta \Bigl(-\int_\Omega \beta \,\nabla_\mathbf{v} f\,d\mathbf{v}\Bigr)
    =
    \int_\Omega \mathbf{v}\,\beta f\,d\mathbf{v}
    - \mathbf{u}_\beta \int_\Omega \beta f\,d\mathbf{v},
\end{equation}
\begin{equation}
    \lambda_\beta \Bigl(-\int_\Omega \beta \mathbf{v}\cdot\nabla_\mathbf{v} f\,d\mathbf{v}\Bigr)
    =
    \int_\Omega \beta\,v^2 f\,d\mathbf{v}
    - \mathbf{u}_\beta \cdot \int_\Omega \mathbf{v}\,\beta f\,d\mathbf{v}.
\end{equation}
To simplify the notation, we introduce the $\beta$-weighted moments:
\begin{equation}
\label{eq:a-p}
    {\boldsymbol{A}}_\beta := -\int_\Omega \beta \,\nabla_\mathbf{v} f\,d\mathbf{v}, 
    \qquad
    \boldsymbol{p}_\beta := \int_\Omega \mathbf{v}\,\beta f\,d\mathbf{v},
\end{equation}
\begin{equation}
\label{eq:n-B-E}
    n_\beta := \int_\Omega \beta f\,d\mathbf{v},
    \qquad
    B_\beta := -\int_\Omega \,\beta \mathbf{v}\cdot\nabla_\mathbf{v} f\,d\mathbf{v},
    \qquad
    E_\beta := \int_\Omega \beta\,v^2 f\,d\mathbf{v}.
\end{equation}
In terms of these moments, the two conditions above can be written as:
\begin{equation}
    \lambda_\beta \,{\boldsymbol{A}}_\beta
    =
    \boldsymbol{p}_\beta - \mathbf{u}_\beta\,n_\beta,
    \qquad
    \lambda_\beta \,B_\beta
    =
    E_\beta - \mathbf{u}_\beta\ \cdot \boldsymbol{p}_\beta.
\end{equation}
Solving for $(\lambda_\beta,\mathbf{u}_\beta)$ gives:
\begin{equation}
    \label{eq:nonlinear_moments}
    \lambda_\beta [f]
    = \frac{n_\beta E_\beta - \bigl(\boldsymbol{p}_\beta\bigr)^2}
           {n_\beta B_\beta - {\boldsymbol{A}}_\beta \cdot \boldsymbol{p}_\beta},
    \qquad
    \mathbf{u}_\beta [f]
    = \frac{1}{n_\beta}\Bigl(\boldsymbol{p}_\beta - \lambda_\beta[f] \,{\boldsymbol{A}}_\beta\Bigr),
\end{equation}
where we have made explicit that the parameters $(\lambda_\beta,\mathbf{u}_\beta)$ are functionals of the distribution function $f$ via the moment integrals. {For stability of the penalization operator, $\lambda_\beta \in (0,\infty)$ as it is an effective temperature. The definition in \eqref{eq:nonlinear_moments} guarantees this property for arbitrary $f$, as shown in \ref{app:lb_moms_well_posedness}.}

In the case of $(\mathbf{u}_\beta^{n},\lambda_\beta^{n})$, $f^n$ is known and the system in \eqref{eq:nonlinear_moments} offers a straightforward method to find them. However, for $(\mathbf{u}_\beta^{n+1},\lambda_\beta^{n+1})$ the situation is very different, because $f^{n+1}$ is not known, as it requires inverting \eqref{eq:penalized_RFP}. It therefore requires a nonlinear iteration where we converge the penalization moments along with the solution of \eqref{eq:penalized_RFP}. We describe the nonlinear solution procedure later in this study.

\subsection{Null-space preservation property of conservative penalization operators $\mathcal{L}_\beta$}

We show next that, as proposed, the penalization operators for the RFP are strictly null-space preserving, i.e., they vanish exactly when $f=f^M$, where $f^M$ is the equilibrium Maxwellian defined in \eqref{Maxwellian_kernel}.
To show this, it is sufficient to show that $\mathbf{u}_\beta \rightarrow \mathbf{u}^M$ and $\lambda_\beta \rightarrow v_{th,M}^2$ as $f\rightarrow f^M$, since:
\begin{equation}
    \nabla_\mathbf{v} f^M + \tfrac{\mathbf{v}-\mathbf{u}^M}{v_{th,M}^2}f^M = \mathbf{0}.
\end{equation}
We begin by noting that the gradient of the Maxwellian is:
\begin{equation}
    \nabla_\mathbf{v} f^M
    =
    -\frac{(\mathbf{v}-\mathbf{u}^M)}{v_{th,M}^2}
    f^M.
\end{equation}
Using this in the definitions of $\bar{\boldsymbol{A}}_\beta$ and $B_\beta$ in the previous section, we readily obtain:
\begin{equation}
    B_\beta^M
    =
    \frac{1}{v_{th,M}^2} \bigl( E_\beta - \mathbf{u}^M \cdot {\mathbf{p}_\beta} \bigr),
    \qquad
    {\boldsymbol{A}}_\beta^M
    = \frac{1}{v_{th,M}^2}  \bigl(
    \boldsymbol{p}_\beta - \mathbf{u}^M \,n_\beta \bigr) .
\end{equation}
Substituting these identities into the expression for \(\lambda_\beta\)
and \(\mathbf{u}_\beta\)
in~\eqref{eq:nonlinear_moments}, we find:
\begin{equation}
    \lambda_\beta = v_{th,M}^2 \,\,;\,\,
    \mathbf{u}_\beta = \mathbf{u}^M.
\end{equation}
i.e., the penalization-operator moments revert back to the Maxwellian
ones, as we sought. This proves strict null-space preservation of the penalized semi-discrete RFP temporal update.
In the particular case of cylindrical coordinates with azimuthal symmetry, the moment unknowns reduce to just $(\lambda_\beta,u_{\parallel,\beta})$.

\section{Structure-preserving spatial discretization of the EIN-penalized RFP equation}
\label{sec:RFP_spatial_discretization}

The previous section proposes a semi-discrete structure-preserving formulation of the penalized RFP equation. The formulation is conservative and equilibrium-preserving. These properties may not survive the velocity-space discretization unless care is taken. In this section, we outline a discretization approach that satisfies all the desired properties in the discrete, namely: maximum principle, equilibrium (null space) preservation, and strict conservation. In what follows, we consider cylindrical geometry and denote $f_{i,j}(t) = f(v_{\perp,i},v_{\parallel,j},t)$, with grids
$v_{\perp,i}=(i-\tfrac{1}{2})\Delta v_\perp$ and
$v_{\parallel,j}=v_{\parallel,\min}+(j-\tfrac{1}{2})\Delta v_\parallel$,
in the perpendicular and parallel directions with respect to some applied electric/magnetic field, respectively.

\subsection{Structure-preserving discretization of the penalization operator $\mathcal{L}_\beta$}

In cylindrical velocity coordinates, it is convenient to write the penalization operator in flux form:
\[
    \mathcal{L}_\beta(f)
    =
    \frac{1}{v_\perp}\,\partial_{v_\perp}\!\bigl(v_\perp J_{\beta,\perp}\bigr)
    + \partial_{v_\parallel} J_{\beta,\parallel},
\]
with penalization flux:
\[
    \mathbf{J}_\beta
    =
    \begin{pmatrix}
      J_{\beta,\perp} \\[0.2em]
      J_{\beta,\parallel}
    \end{pmatrix}
    =
    \beta(v_\perp,v_\parallel)
    \Bigl(
      \nabla_\mathbf{v} f
      + \tfrac{\mathbf{v}-\mathbf{u}_\beta}{\lambda_\beta}f
    \Bigr),
\]
We denote $\beta_{i,j} = \beta(v_{\perp,i},v_{\parallel,j})$ and define
face-centered values by simple averaging,
\[
    \beta_{i\pm 1/2,j}
    = \tfrac{1}{2}\bigl(\beta_{i,j}+\beta_{i\pm1,j}\bigr),
    \qquad
    \beta_{i,j\pm 1/2}
    = \tfrac{1}{2}\bigl(\beta_{i,j}+\beta_{i,j\pm1}\bigr).
\]
A mass-conservative finite-difference discretization of $\mathcal{L}_\beta$
at $(i,j)$ then reads:
\begin{equation}
    (\mathcal{L}_\beta f)_{i,j}
    =
    \frac{1}{v_{\perp,i}\,\Delta v_\perp}
    \bigl(
        v_{\perp,i+1/2} J_{\beta,\perp,i+1/2,j}
        - v_{\perp,i-1/2} J_{\beta,\perp,i-1/2,j}
    \bigr)
    +
    \frac{1}{\Delta v_\parallel}
    \bigl(
        J_{\beta,\parallel,i,j+1/2}
        - J_{\beta,\parallel,i,j-1/2}
    \bigr),
\end{equation}
where the fluxes at faces are split into diffusive and drift contributions,
\[
    J_{\beta,\perp,i+1/2,j}
    =
    \beta_{i+1/2,j}
    \frac{f_{i+1,j}-f_{i,j}}{\Delta v_\perp}
    +
    a_{\beta,\perp,i+1/2,j}\,f_{i+1/2,j},
    \qquad
    a_{\beta,\perp,i+1/2,j}
    =
    \beta_{i+1/2,j}\,
    \frac{v_{\perp,i+1/2}}{\lambda_\beta},
\]
\begin{equation}
\label{eq:j_par_beta_def}
    J_{\beta,\parallel,i,j+1/2}
    =
    \beta_{i,j+1/2}
    \frac{f_{i,j+1}-f_{i,j}}{\Delta v_\parallel}
    +
    a_{\beta,\parallel,i,j+1/2}\,f_{i,j+1/2},
    \qquad
    a_{\beta,\parallel,i,j+1/2}
    =
    \beta_{i,j+1/2}\,
    \frac{v_{\parallel,j+1/2}-u_{\parallel,\beta}}{\lambda_\beta}.
\end{equation}
Here, $v_{\perp,i\pm 1/2} = v_{\perp,i} \pm \tfrac{1}{2}\Delta v_\perp$,
$v_{\parallel,j\pm 1/2} = v_{\parallel,j} \pm \tfrac{1}{2}\Delta v_\parallel$, and $f_{i+1/2,j}$ and $f_{i,j+1/2}$ denote the unknown values of $f$ at
the faces, which we approximate using the \emph{classical}
Chang-Cooper discretization \cite{chang1970practical} in each coordinate direction, given by:
\begin{align}
\label{solution_at_faces_CC_1}
    f_{i+1/2,j}
    &= (1-\delta_{\perp,i+1/2,j})\,f_{i+1,j}
       + \delta_{\perp,i+1/2,j}\,f_{i,j}, \\
       \label{solution_at_faces_CC_2}
    \delta_{\perp,i+1/2,j}
    &= \frac{1}{\Delta v_\perp w_{\perp,i+1/2,j}}
       - \frac{1}{\exp\!\bigl(\Delta v_\perp w_{\perp,i+1/2,j}\bigr)-1}, 
    \qquad
    w_{\perp,i+1/2,j}
    = \frac{v_{\perp,i+1/2}}{\lambda_\beta},
    \\[0.4em]
    \label{solution_at_faces_CC_3}
    f_{i,j+1/2}
    &= (1-\delta_{\parallel,i,j+1/2})\,f_{i,j+1}
       + \delta_{\parallel,i,j+1/2}\,f_{i,j}, \\
       \label{solution_at_faces_CC_4}
    \delta_{\parallel,i,j+1/2}
    &= \frac{1}{\Delta v_\parallel w_{\parallel,i,j+1/2}}
       - \frac{1}{\exp\!\bigl(\Delta v_\parallel w_{\parallel,i,j+1/2}\bigr)-1}, 
    \qquad
    w_{\parallel,i,j+1/2}
    = \frac{v_{\parallel,j+1/2}-u_{\parallel,\beta}}{\lambda_\beta}.
\end{align}

The Chang-Cooper discretization yields a maximum-principle- and analytical-null-space-preserving discretization of the penalization operator (when $\mathbf{u}_\beta$ and $\lambda_\beta$ are defined as in Sec. \ref{sec:cons_penalization}). 
However, the discretization is not necessarily conservative yet. For strict discrete conservation, we also need the moment integrals for $\mathbf{u}_\beta$ and $\lambda_\beta$ in \eqref{eq:lb_cons_constraints} to be discretized correctly. We derive a suitable discrete form of these integrals next.

We begin from the first-moment integral, which in cylindrical geometry simplifies to considering just the parallel component (the perpendicular component vanishes in the integral due to azimuthal symmetry).
Discretizing the $v_\parallel$ integral and performing a discrete integration by parts yields:
\begin{align}
    \int_\Omega v_\parallel \mathcal{L}_\beta f d\mathbf{v}
    &\approx 2\pi\,\Delta v_\perp \Delta v_\parallel
    \sum_{i,j}
    \Biggl[
        v_{\perp,i} v_{\parallel,j}\,
        \frac{
            \beta_{i,j+1/2} J_{\beta,\parallel,i,j+1/2}
            - \beta_{i,j-1/2} J_{\beta,\parallel,i,j-1/2}
        }{\Delta v_\parallel}
    \Biggr] \\
    &= 2\pi\,\Delta v_\perp \Delta v_\parallel
    \sum_{i,j}
    v_{\perp,i}\,\beta_{i,j+1/2}\,J_{\beta,\parallel,i,j+1/2}\,
    \frac{v_{\parallel,j}-v_{\parallel,j+1}}{\Delta v_\parallel} \\
    &= -2\pi\,\Delta v_\perp \Delta v_\parallel
    \sum_{i,j}
    v_{\perp,i}\,\beta_{i,j+1/2}\,J_{\beta,\parallel,i,j+1/2}=0.
\label{eq:lb_par_int}
\end{align}
Introducing the discrete definition of $J_{\beta,\parallel}$ in \eqref{eq:j_par_beta_def} into \eqref{eq:lb_par_int} yields the discrete definitions of
$A_{\beta,\parallel}$, $p_{\beta,\parallel}$, and $n_\beta$ sought:
\begin{align}
\label{A_beta_moment}
    A_{\beta,\parallel}
    &= -2\pi\,\Delta v_\perp \Delta v_\parallel
       \sum_{i,j} v_{\perp,i}\,\beta_{i,j+1/2}\,
       \frac{f_{i,j+1}-f_{i,j}}{\Delta v_\parallel}, \\[0.3em]
       \label{p_beta_moment}
    p_{\beta,\parallel}
    &= 2\pi\,\Delta v_\perp \Delta v_\parallel
       \sum_{i,j} v_{\perp,i}\,\beta_{i,j+1/2}\,v_{\parallel,j+1/2}\,
       f_{i,j+1/2}, \\[0.3em]
       \label{n_beta_moment}
    n_\beta
    &= 2\pi\,\Delta v_\perp \Delta v_\parallel
       \sum_{i,j} v_{\perp,i}\,\beta_{i,j+1/2}\,f_{i,j+1/2}.
\end{align}
Discretizing the second-moment identity
\[
    \int_\Omega v^2\,
    \mathcal{L}_\beta f\,d\mathbf{v}
    = 0,
\]
and manipulating it similarly by performing a numerical integration by parts, we obtain the following discrete definitions of
$B_\beta$ and $E_\beta$:
\begin{align}
\label{B_beta_moment}
    B_\beta
    &= -2\pi\,
       \sum_{i,j} \Biggl[
           \Delta v_\parallel v_{\perp,i+1/2}^2\,\beta_{i+1/2,j}\,
           (f_{i+1,j}-f_{i,j})
           \;+\;
           \Delta v_\perp v_{\perp,i}v_{\parallel,j+1/2}\,\beta_{i,j+1/2}\,
           (f_{i,j+1}-f_{i,j})
       \Biggr], \\[0.3em]
       \label{E_beta_moment}
    E_\beta
    &= 2\pi\,\Delta v_\perp \Delta v_\parallel
       \sum_{i,j} \Bigl[
           v_{\perp,i+1/2}^3\,\beta_{i+1/2,j}\,
           f_{i+1/2,j}
           \;+\;
           v_{\perp,i}v_{\parallel,j+1/2}^2\,\beta_{i,j+1/2}\,
           f_{i,j+1/2}
       \Bigr].
\end{align}
With these discrete definitions, one can find the corresponding $\lambda_\beta$ and $u_{\beta,\parallel}$ parameters according to \eqref{eq:nonlinear_moments} such that the penalization operator is strictly conserving of momentum and energy. Note that the discrete solution values at the cell faces $f_{i+1/2,j}$ and $f_{i,j+1/2}$ are obtained using the classical Chang-Cooper interpolation, \eqref{solution_at_faces_CC_1}-\eqref{solution_at_faces_CC_4}, which has a maximum principle and is equilibrium (null-space) preserving. However, because the $\delta$-parameter in the Chang-Cooper discretization itself depends on $\lambda_\beta$ and $u_{\parallel,\beta}$ in a nonlinear fashion (through $\omega_\perp$ and $\omega_\parallel$), the discrete version of \eqref{eq:nonlinear_moments} is nonlinear in $u_{\parallel,\beta}$ and $\lambda_\beta$ and requires iteration, regardless of whether $f_{i,j}$ is known on the mesh or not. We will discuss this in more detail later in this study. 

\subsection{Structure-preserving spatial discretization of the RFP operator}

We consider next the structure-preserving discretization of the RFP collision operator, also in
cylindrical velocity coordinates $(v_\perp,v_\parallel)$,
\begin{equation}
\label{eq:RFP_discretization}
   \mathcal{C}(f)
    \;=\;
    \nabla_\mathbf{v}\!\cdot\mathbf{J}
    \;=\;
    \frac{1}{v_\perp}\,\partial_{v_\perp}\!\bigl(v_\perp J_\perp\bigr)
    + \partial_{v_\parallel} J_\parallel,
\end{equation}
with flux $\mathbf{J}= (\gamma \bar{\bar{I}} + \bar{\bar{\epsilon}}){\cal D} \cdot\nabla_\mathbf{v} f -\mathbf{A}f =(J_\perp,J_\parallel)^\top$ and a
symmetric positive-definite diffusion tensor
\[
    (\gamma \bar{\bar{I}} + \bar{\bar{\epsilon}}){\cal D}(v_\perp,v_\parallel)
    =
    \begin{pmatrix}
        D_{\perp\perp} & D_{\perp\parallel}\\[0.2em]
        D_{\parallel\perp} & D_{\parallel\parallel}
    \end{pmatrix},
    \qquad
    D_{\perp\parallel}=D_{\parallel\perp}.
\]
Here:
\[
\bar{\bar{\epsilon}}
=
\begin{pmatrix}
\epsilon_\parallel & 0 \\
0                  & 0
\end{pmatrix}.
\]
The parameters $\gamma$ and $\epsilon_\parallel$ ensure exact discrete conservation of mass, momentum, and energy, as proposed in \cite{taitano2015mass}, and are defined discretely later in this study. 

For strict maximum-principle and null-space (equilibrium) preservation, we will combine the adaptive timestepping and advectionalization approaches introduced earlier with a novel adaption of the Chang-Cooper discretization, discussed below. The resulting discretization will satisfy a maximum principle, conservation, and preservation of the analytical null-space strictly.
We begin with a mass-conservative finite-difference discretization of the cylindrical divergence
operator, given by:
\begin{equation}
  \mathcal{C}(f)_{i,j}
    =
    \frac{1}{v_{\perp,i}\,\Delta v_\perp}
        \bigl(
            v_{\perp,i+1/2} J_{\perp,i+1/2,j}
            - v_{\perp,i-1/2} J_{\perp,i-1/2,j}
        \bigr)
    +
    \frac{1}{\Delta v_\parallel}
        \bigl(
            J_{\parallel,i,j+1/2}
            - J_{\parallel,i,j-1/2}
        \bigr).
\end{equation}
As we did before for the anisotropic heat equation, we split the fluxes into a diagonal part and an
off-diagonal part:
\[
    J_\perp = J_\perp^{\mathrm{diag}} + J_\perp^{\mathrm{aniso}},
    \qquad
    J_\parallel = J_\parallel^{\mathrm{diag}} + J_\parallel^{\mathrm{aniso}}.
\]
The diagonal part is discretized with centered finite differences, e.g.:
\begin{equation}
\label{fluxes_RFP_diag}
    J_{\perp,i+1/2,j}^{\mathrm{diag}}
    = D_{\perp\perp,i+1/2,j}\,
      \frac{f_{i+1,j}-f_{i,j}}{\Delta v_\perp},
    \qquad
    J_{\parallel,i,j+1/2}^{\mathrm{diag}}
    = D_{\parallel\parallel,i,j+1/2}\,
      \frac{f_{i,j+1}-f_{i,j}}{\Delta v_\parallel},
\end{equation}
with suitable averages for $D_{\perp\perp}$ and $D_{\parallel\parallel}$.
The off-diagonal anisotropic contributions are advectionalized as:
\begin{equation}
\label{fluxes_RFP_aniso}
    J_{\perp,i+1/2,j}^{\mathrm{aniso}} = a_{\perp,i+1/2,j}\,f_{i+1/2,j},
    \qquad
    J_{\parallel,i,j+1/2}^{\mathrm{aniso}} = a_{\parallel,i,j+1/2}\,f_{i,j+1/2},
\end{equation}
where as before the effective velocities are defined from the derivatives of
$\ln(f)$ and the friction velocity at faces as:
\begin{align}
    a_{\perp,i+1/2,j}
    &=
    \bigl(D_{\perp\parallel}\partial_{v_\parallel}\ln f - A_\perp \bigr)_{i+1/2,j}
    =
    D_{\perp\parallel,i+1/2,j}\,
    \frac{\ln f_{i+1/2,j+1/2} - \ln f_{i+1/2,j-1/2}}{\Delta v_\parallel}
    -A_{\perp,{i+1/2,j}} ,
    \\[0.3em]
    a_{\parallel,i,j+1/2}
    &=
    \bigl(D_{\perp\parallel}\partial_{v_\perp}\ln f- A_\parallel \bigr)_{i,j+1/2}
    =
   D_{\perp\parallel,i,j+1/2}\,
   \frac{\ln f_{i+1/2,j+1/2} - \ln f_{i-1/2,j+1/2}}{\Delta v_\perp}
   -A_{\parallel,{i,j+1/2}}.
   \label{coeffs_RFP_aniso}
\end{align}
Here, the advection and diffusion coefficients are computed from two
Poisson potential solves using a second-order centered-difference
discretization with suitable far-field boundary conditions \cite{taitano2016adaptive}. More
details are provided in~\ref{Appendix:Poisson_Potentials}. Corner values such as $f_{i+1/2,j+1/2}$ are approximated by the
average of the four surrounding cell centers, e.g.
\[
    f_{i+1/2,j+1/2}
    \approx \frac{1}{4}\bigl(
        f_{i,j} + f_{i+1,j} + f_{i,j+1} + f_{i+1,j+1}
    \bigr),
\]
and similarly for $f_{i+1/2,j-1/2}$ and $f_{i-1/2,j+1/2}$. The face values
$f_{i+1/2,j}$ and $f_{i,j+1/2}$ needed for the advective fluxes are reconstructed from cell-centered values $\{f_{i,j}\}$ by a modified Chang-Cooper interpolation that guarantees both a maximum principle and analytical equilibrium preservation. We discuss this next.

\subsubsection{Maximum-principle- and null-space-preserving RFP spatial discretization}\label{MP_discretization_RFP}

In each coordinate direction, the classical Chang-Cooper scheme interpolates
the solution from cell centers to the cell faces for the advective flux in \eqref{fluxes_RFP_aniso} according to:
\begin{align}
    f_{i+1/2,j}
    &= (1-\delta_{\perp,i+1/2,j})\,f_{i+1,j}
       + \delta_{\perp,i+1/2,j}\,f_{i,j}, \\
    \delta_{\perp,i+1/2,j}
    &= \frac{1}{\Delta v_\perp w_{\perp,i+1/2,j}}
       - \frac{1}{\exp\!\bigl(\Delta v_\perp w_{\perp,i+1/2,j}\bigr)-1}, 
    \qquad
    w_{\perp,i+1/2,j}
    = \frac{a_{\perp,i+1/2,j}}{D_{\perp\perp,i+1/2,j}},
    \\[0.4em]
    f_{i,j+1/2}
    &= (1-\delta_{\parallel,i,j+1/2})\,f_{i,j+1}
       + \delta_{\parallel,i,j+1/2}\,f_{i,j}, \\
    \delta_{\parallel,i,j+1/2}
    &= \frac{1}{\Delta v_\parallel w_{\parallel,i,j+1/2}}
       - \frac{1}{\exp\!\bigl(\Delta v_\parallel w_{\parallel,i,j+1/2}\bigr)-1}, 
    \qquad
    w_{\parallel,i,j+1/2}
    = \frac{a_{\parallel,i,j+1/2}}{D_{\parallel\parallel,i,j+1/2}}.
\end{align}
For the Maxwellian steady state, the Chang-Cooper weights $(\omega_\perp,\omega_\parallel)$ can be found analytically, either from the analytical solutions of the Rosenbluth potentials, or (more straightforwardly) from the RFP equilibrium condition:
\begin{equation}
\label{eq:tensor_in_equilbirium}
    \mathcal{D}^M\cdot\nabla f^M - \mathbf{A}^M f^M = 0
    \;\;\Longrightarrow\;\;
    -\mathcal{D}^M\cdot\,\frac{\mathbf{v}-\mathbf{u}^M}{v_{th,M}^2} = \mathbf{A}^M,
\end{equation}
since $\nabla f^M / f^M = -(\mathbf{v}-\mathbf{u}^M)/v_{th,M}^2$ for a Maxwellian $f^M$ with thermal velocity $v_{th,M}$ and drift $\mathbf{u}^M$.
In cylindrical coordinates $(v_\perp,v_\parallel)$, this relation reads:
\begin{equation}
    -\frac{1}{v_{th,M}^2}
    \begin{bmatrix}
        D^M_{\perp\perp}\,v_\perp + D^M_{\perp\parallel}\,(v_\parallel - u^M_\parallel) \\[0.3em]
        D^M_{\perp\parallel}\,v_\perp + D^M_{\parallel\parallel}\,(v_\parallel - u^M_\parallel)
    \end{bmatrix}
    =
    \begin{bmatrix}
        A^M_\perp \\[0.3em]
        A^M_\parallel
    \end{bmatrix}.
\end{equation}
Using this result, one obtains the following analytical expressions for the Chang-Cooper weights for $f^M$:
\begin{align}
    w^M_\perp
    &= \frac{a_\perp}{D_{\perp\perp}}
     = \frac{D_{\perp\parallel}\,\partial_\parallel \ln(f^M) - A_\perp}{D_{\perp\perp}} \\[0.3em]
    &= \frac{
        -D_{\perp\parallel}\,\dfrac{v_\parallel - u_\parallel^M}{v_{th,M}^2}
        + D_{\perp\perp}\,\dfrac{v_\perp}{v_{th,M}^2}
        + D_{\perp\parallel}\,\dfrac{v_\parallel - u_\parallel^M}{v_{th,M}^2}
    }{D_{\perp\perp}} \\[0.3em]
    &= \frac{v_\perp}{v_{th,M}^2},
\end{align}
and:
\begin{equation}
    w^M_\parallel = \frac{v_\parallel - u_\parallel^M}{v_{th,M}^2}.
\end{equation}

However, in our discrete setting, the advection and diffusion coefficients
$a_{\perp/\parallel}$ and $D$ are computed numerically from the Rosenbluth
potentials, and the resulting discrete weights $w_\perp = a_\perp/D_{\perp\perp}$,
$w_\parallel = a_\parallel/D_{\parallel\parallel}$, will not exactly satisfy
these equilibrium identities. This discrepancy will prevent the strict preservation of the
Maxwellian steady state.
To mitigate this, we modify the Chang-Cooper weights for the RFP equation by combining the
numerically computed ${w_\perp}$ and $w_\parallel$ with the analytically defined equilibrium
weights ${w}^M_\perp$ and  ${w}^M_\parallel$ as (see \ref{app:rfp-chang-cooper}):
\begin{equation}
\label{modified_CC}
    {\theta_{\perp}}
    =
    \frac{1}{{\Delta v_\perp w_\perp}}
    - \frac{1}{\exp\!\bigl(\Delta v_\perp {w}^M_\perp\bigr)-1}, \quad   {\theta_{\parallel}}
    =
    \frac{1}{{\Delta v_\parallel w_\parallel}}
    - \frac{1}{\exp\!\bigl(\Delta v_\parallel {w}^M_\parallel\bigr)-1}.
\end{equation}
However, as defined, the proposed RFP Chang-Cooper weights ${\theta}$ may no longer remain bounded in the interval $[0,1]$, as is the
case for the classical scheme (and as needed for stability). To enforce this, we propose to only employ this
modified interpolation when the numerical solution is already very close to the
equilibrium Maxwellian, as measured by the condition:
\[
    \bigl\|
        \mathbf{w}
        - {\mathbf{w}^M}
    \bigr\|
    < \varepsilon_{CC},
\]
for some user-defined tolerance $\varepsilon_{CC} > 0$. Otherwise, the standard Chang-Cooper interpolation is used. As an additional safety net, we ensure that the modified Chang-Cooper weight remains bounded in $[0,1]$ at all times by enforcing:
\begin{equation}
  \theta_{\perp,\parallel} \;\leftarrow\;
  \begin{cases}
    \displaystyle \max\bigl(0,\min(\theta,\tfrac{1}{2})\bigr),
      & \text{if } w_{\perp,\parallel} > 0,\\[0.6em]
    \displaystyle \max\bigl(\tfrac{1}{2},\min(\theta,1)\bigr),
      & \text{if } w_{\perp,\parallel} < 0,
  \end{cases}
\end{equation}
{To ensure that the established equilibrium properties are preserved at the boundaries, we fill ghost cells using the analytical Maxwellian \( f^M \).}

\subsubsection{Strict discrete conservation strategy in RFP}

For the RFP operator, we employ a strictly conservative discretization
following~\cite{taitano2015mass}. We discretize
$\mathbf{J}^{\mathrm{diag}} = [J_\perp^{\mathrm{diag}}, J_\parallel^{\mathrm{diag}}]$
according to~\eqref{fluxes_RFP_diag}, and
$J^{\mathrm{aniso}} = [J_\perp^{\mathrm{aniso}}, J_\parallel^{\mathrm{aniso}}]$
according to~\eqref{fluxes_RFP_aniso}-\eqref{coeffs_RFP_aniso}.  
With these fluxes, we compute the conservation moments:
\[
\gamma
=
\frac{
\left\langle \mathbf{v}, \mathbf{J}^{\mathrm{aniso}} \right\rangle_{{\mathbf{v}}}
-
\epsilon_{\parallel}
\left\langle {v}_\parallel, J^{\mathrm{diag}}_{\parallel} \right\rangle_{u_{\parallel}}^{+\infty}
}{
\left\langle \mathbf{v}, \mathbf{J}^{\mathrm{diag}} \right\rangle_{{\mathbf{v}}}},
\]
and
\[
\epsilon_{\parallel}
=
\begin{cases}
\displaystyle
\frac{
\left\langle 1, J^{\mathrm{aniso}}_{\parallel} \right\rangle_{{\mathbf{v}}}
-
\gamma \left\langle 1, J^{\mathrm{diag}}_{\parallel} \right\rangle_{{\mathbf{v}}}
}{
\left\langle 1, J^{\mathrm{diag}}_{\parallel} \right\rangle_{u_{\parallel}}^{+\infty}
},
& \text{if } v_{\parallel} \ge u_{\parallel},\\[1ex]
0, & \text{otherwise}.
\end{cases}
\]
We then assemble the RFP collision operator as:
\begin{equation}
  \mathcal{C}(f)_{i,j}
    =
    \frac{\gamma}{v_{\perp,i}\,\Delta v_\perp}
        \bigl(
            v_{\perp,i+1/2} J_{\perp,i+1/2,j}
            - v_{\perp,i-1/2} J_{\perp,i-1/2,j}
        \bigr)
    +
    \frac{\epsilon_\parallel + \gamma}{\Delta v_\parallel}
        \bigl(
            J_{\parallel,i,j+1/2}
            - J_{\parallel,i,j-1/2}
        \bigr).
\end{equation}
At equilibrium, the modified Chang-Cooper discretization preserves the
steady state analytically, so the numerical parameters satisfy
$\epsilon_\parallel = 0$ and $\gamma = 1$ numerically to machine precision.

\subsubsection{Nonlinear solution strategy for the penalized RFP discretized system}

At each time-step, we first determine the conservative penalization parameters, $(\mathbf{u}_\beta^{n},\lambda_\beta^{n})[f^{n}]$, using Algorithm \ref{alg:conservative-penalized-update-RHS}. Since $f^n$ is available, this reduces to solving a two-equation nonlinear system, with the nonlinearity resulting from the dependence of the Chang-Cooper discretization of the penalization-operator fluxes on $(\mathbf{u}_\beta^{n},\lambda_\beta^{n})$. 

\begin{algorithm}[t]
  \caption{Solution algorithm for ($\lambda_\beta^n,u_{\parallel,\beta}^n$)}
  \label{alg:conservative-penalized-update-RHS}
  \begin{algorithmic}[1]
    \Require $f^n$
    \State Choose an initial guess, e.g.
    \[
        \Xi_{\beta,(0)}^{n}=({u}_\parallel,v_{th,M}^2).
    \]
    \For{$k = 0,1,2,\dots$ until convergence}
        \State {Compute Chang-Cooper weights from from $\Xi_{\beta,(k)}^{n}$, and interpolate solution at faces.}
        \State Update the $\beta$-weighted moments $\Xi_{\beta,(k+1)}^{n}$ with $f^n$
        via~\eqref{A_beta_moment}-\eqref{E_beta_moment}.% 
        \State Accelerate with AA:
        \[
            \Xi_{\beta,(k+1)}^{n}
            =
            \Xi_{\beta,(k)}^{n}
            +
            AA\bigl(\Xi_{\beta,(k+1)}^n-\Xi_{\beta,(k)}^{n}\bigr).
        \]
        \If{
        $
            \bigl\|\mathbf{u}_{\beta,(k+1)}^{\,n}
                - \mathbf{u}_{\beta,(k)}^{\,n}\bigr\|
            < \varepsilon_{AA},
            \qquad
            \bigl|\lambda_{\beta,(k+1)}^{\,n}
                - \lambda_{\beta,(k)}^{\,n}\bigr|
            < \varepsilon_{AA}
        $}
            \State \textbf{stop}.
        \EndIf
    \EndFor
  \end{algorithmic}
\end{algorithm}

Using these parameters, our solution strategy for the penalized RFP equation~\eqref{eq:penalized_RFP} is outlined in Algorithm \ref{alg:conservative-penalized-update}. This algorithm iteratively solves the nonlinear problem for $(\mathbf{u}_\beta^{n+1},\lambda_\beta^{n+1})[f^{n+1}]$. 
This second system is more involved, as $f^{n+1}$ is not known, and needs to be updated during the iteration. This is accomplished as described in the Algorithm. The iteration is a Picard algorithm, accelerated with Anderson Acceleration (AA) \cite{anderson1965iterative,walker2011anderson} for faster performance. In our implementation, we use five AA states as our acceleration subspace. In our numerical experiments, AA is able to improve performance over Picard by a factor of two to four, depending on timestep size, resulting in fewer than a handful of iterations for relatively tight tolerances, as documented in the numerical section.

\begin{algorithm}[t]
  \caption{Solution algorithm for ($\lambda_\beta^{n+1},u_{\parallel,\beta}^{n+1},f^{n+1}$)}
  \label{alg:conservative-penalized-update}
  \begin{algorithmic}[1]
    \Require $f^n$
    \State Using Algorithm \ref{alg:conservative-penalized-update-RHS}, solve for 
    \[
        \Xi_\beta^n=(\mathbf{u}_\beta^{n},\lambda_\beta^{n}).
    \]
    \State Choose an initial guess for the next-step parameters, e.g.
    \[
        \Xi_{\beta,(0)}^{n+1}=\Xi_\beta^n.
    \]
    \For{$k = 0,1,2,\dots$ until convergence}
        \State PDF (linear) update: solve \eqref{eq:penalized_RFP} with the proposed SP discretization for $f^\ast$ using
        \[
            \Xi_{\beta,(k)}^{n+1}
            =
            (\mathbf{u}_{\beta,(k)}^{\,n+1},\lambda_{\beta,(k)}^{\,n+1}).
        \]
        \State Update the $\beta$-weighted moments $\Xi_{\beta,(k+1)}^{*}$ with $f^\ast$
        via~\eqref{eq:nonlinear_moments}.
        \State Accelerate with AA:
        \[
            \Xi_{\beta,(k+1)}^{n+1}
            =
            \Xi_{\beta,(k)}^{n+1}
            +
            AA\bigl(\Xi_{\beta,(k+1)}^*-\Xi_{\beta,(k)}^{n+1}\bigr).
        \]
        \If{
        $
            \bigl\|\mathbf{u}_{\beta,(k+1)}^{\,n+1}
                - \mathbf{u}_{\beta,(k)}^{\,n+1}\bigr\|
            < \varepsilon_{AA},
            \qquad
            \bigl|\lambda_{\beta,(k+1)}^{\,n+1}
                - \lambda_{\beta,(k)}^{\,n+1}\bigr|
            < \varepsilon_{AA}
        $}
            \State Set $f^{n+1}=f^\ast$ and \textbf{stop}.
        \EndIf
    \EndFor
  \end{algorithmic}
\end{algorithm}

\section{Numerical results}
\label{sec:numerical_results}
We consider tests of varying complexity, starting from linear anisotropic-diffusion tests, continuing with a linear cross-species collisions with disparate masses (``pitch-angle scattering''), and culminating with several fully nonlinear single-species RFP tests. For all these tests, unless otherwise stated, we {impose zero-flux (homogeneous Neumann) boundary conditions} and set $\varepsilon = \min\!\left(0.05, \frac{\lambda_1-\lambda_2}{4\beta}\right)$ for the adaptive time-stepping scheme, $\epsilon_{CC}=0.1$ for the modified Chang-Cooper scheme, and $\epsilon_{AA}=10^{-10}$ for the Anderson-acceleration tolerances. In order to assess the magnitude of the time step, we introduce the CFL-count parameter
\begin{equation}
    N_{\mathrm{CFL}} = \frac{\Delta t}{\Delta t_{\mathrm{CFL}}}.
\end{equation}
This quantity serves as a metric to gauge how large the chosen time step is compared to the explicitly constrained one. {The penalization parameter $\beta$ is taken as $\beta(\mathbf{x})=\tfrac{1}{2}\lambda_1(\mathbf{x})$, where $\lambda_1(\mathbf{x})$ denotes the largest eigenvalue of the diffusion tensor $D$ computed numerically over the mesh. For the heat equation tests and the linearized RFP model, $\beta$ is computed once at the beginning, since the diffusion tensor is time independent. For the nonlinear RFP tests, $\beta$ is updated at each time step, since the diffusion tensor depends on the evolving solution.}

\subsection{Anisotropic Diffusion tests}
We consider both a constant-diffusion-tensor and a variable-diffusion-tensor problem. For the following test, we use a Cartesian domain $[-1,1]\times[-1,1]$. We fix the grid size to $N_x=N_y=100$. 

\subsubsection{Constant anisotropic diffusion}
We first consider a constant-coefficient anisotropic diffusion tensor. 
The initial condition is a Gaussian bump centered at the origin,
\[
    f_0(\mathbf{x})
    =
    A \exp\!\Bigl(
        -\frac{x^2+y^2}{\sigma^2}
    \Bigr),
    \qquad
    \sigma = 0.1,
    \quad
    A = \frac{1}{\pi \sigma^2}.
\]
The diffusion tensor is constant in space and given by:
\[
    D
    =
    R(\theta)
    \begin{pmatrix}
        \lambda_1 & 0\\[0.2em]
        0 & \lambda_2
    \end{pmatrix}
    R(\theta)^\top,
    \qquad
    \lambda_1 = 1,\quad
    \lambda_2 = 10^{-3},\quad
    \theta = \tfrac{3\pi}{8},
\]
where $R(\theta)$ is the rotation matrix of angle $\theta$. This setup yields a strongly anisotropic diffusion with anisotropy ratio
$\lambda_1/\lambda_2 = 10^3$, oriented at an angle $\theta=3\pi/8$ with respect to
the $x$-axis, and is used to assess the ability of the scheme to resolve
strong, rotated anisotropy while preserving positivity and accuracy. 
The Fourier analysis used to design the time-adaptive scheme, in section \ref{Sec:Time_Adaptivity}, is exact for this test, since the analysis assumes a constant diffusion tensor.
\begin{figure}
    \centering    \includegraphics[width=0.8\linewidth]{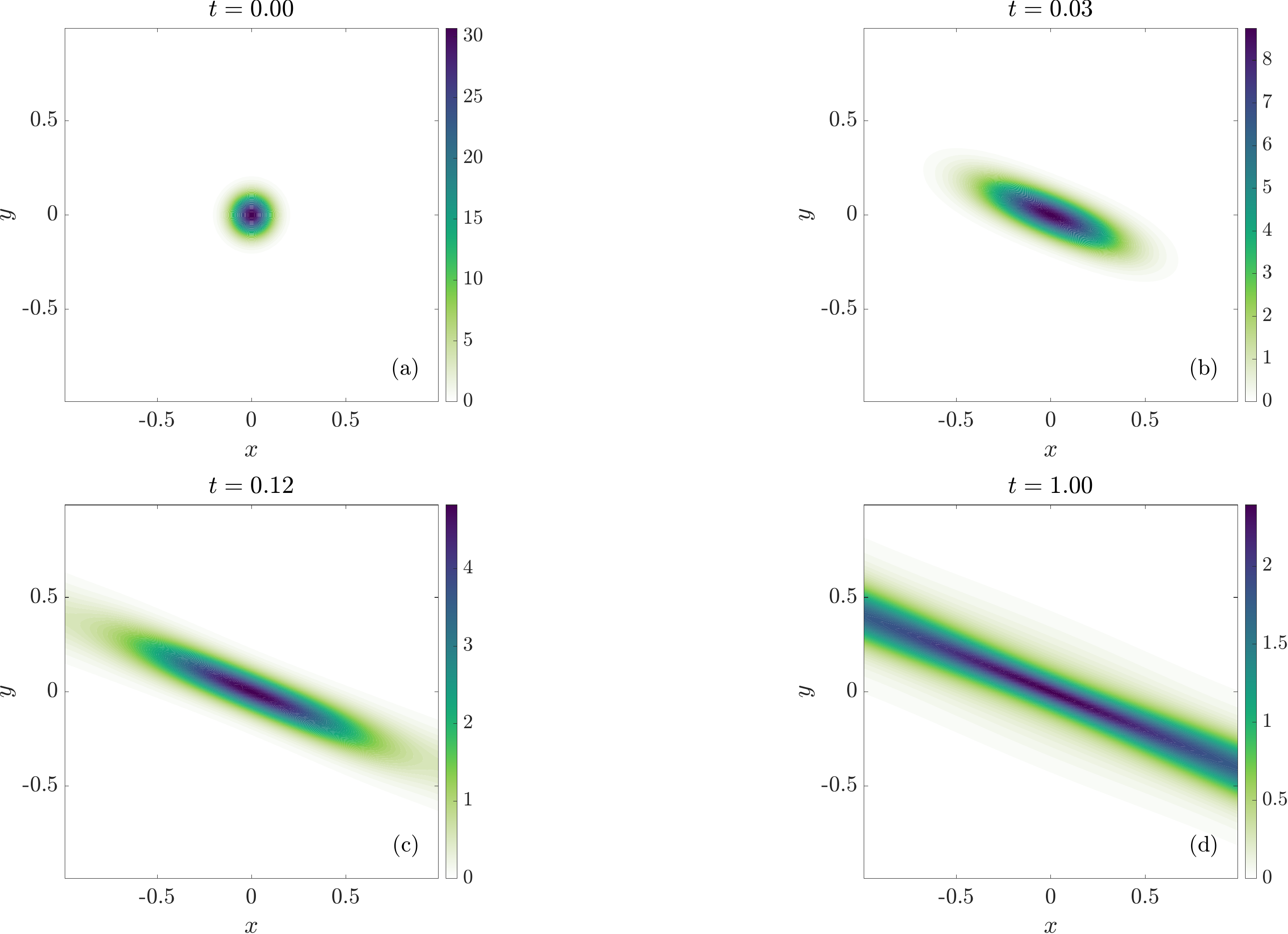}
    \caption{Evolution of the solution for the anisotropic ring test at four representative times.}
    \label{fig:band_sim}
\end{figure}

Figure~\ref{fig:band_sim} shows the temporal evolution of the solution. The initial Gaussian profile diffuses along the rotated principal axes specified by $\theta$, progressively forming a stretched band aligned with these anisotropic directions.

\begin{figure}
    \centering
    \includegraphics[width=\linewidth]{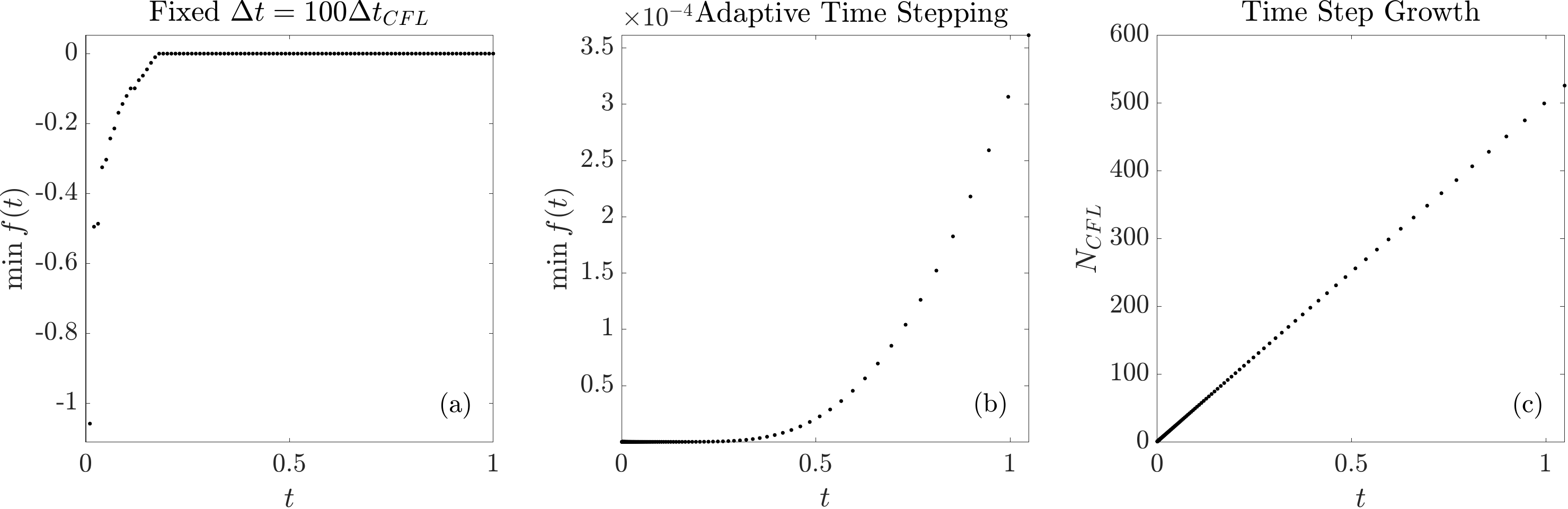}
   \caption{(a) Minimum value of the solution in time for a fixed time-step ratio $\Delta t = 100\,\Delta t_{\mathrm{CFL}}$. (b) Minimum value of the solution in time for the adaptive time-stepping strategy. (c) Evolution of the adaptive time-step ratio $N_{CFL}=\Delta t_n / \Delta t_{\mathrm{CFL}}$ as a function of time.}

    \label{fig:band_negativity}
\end{figure}

Figure~\ref{fig:band_negativity} (a) displays the evolution of the adaptive timestep as a function of the time-step index. The timestep grows  exponentially as the iteration proceeds, as expected from our analysis. Figures~\ref{fig:band_negativity} (b)-(c) compare the smallest value of the solution for a fixed time step with $N_{CFL}=100$ and for the adaptive time-stepping strategy, respectively. In the fixed time-step case, the solution exhibits pronounced negativity at early times, due to the contribution from the anisotropic kernel for unresolved modes. In contrast, the adaptive time-stepping simulation remains positive as the time step is increased progressively up to $N_{CFL}=500$, demonstrating the effectiveness of the prescription.

\subsubsection{Variable anisotropic diffusion (ring test)}
We consider next a variable anisotropic diffusion tensor that generates purely rotational diffusion \cite{sharma2007preserving,deka2022exponential,crouseilles2015comparison}. The initial condition is taken as a normalized Gaussian bump centered at $(-0.6,0)$,
\[
    f_0(\mathbf{x})
    =
    A \exp\!\Bigl(
        -\frac{(x+0.6)^2}{\sigma^2}
        -\frac{y^2}{\sigma^2}
    \Bigr),
    \qquad
    \sigma = 0.1,
    \quad
    A = \frac{1}{\pi \sigma^2}.
\]
The anisotropic diffusion tensor is given by:
\[
    D(\mathbf{x})
    =
    \begin{pmatrix}
        y^2   & -x y\\[0.2em]
       -x y  & x^2
    \end{pmatrix}.
\]
This tensor has eigenvalues $0$ and $x^2+y^2$. The diffusion strength increases quadratically with the distance from the origin.
\begin{figure}
    \centering
    \includegraphics[width=0.7\linewidth]{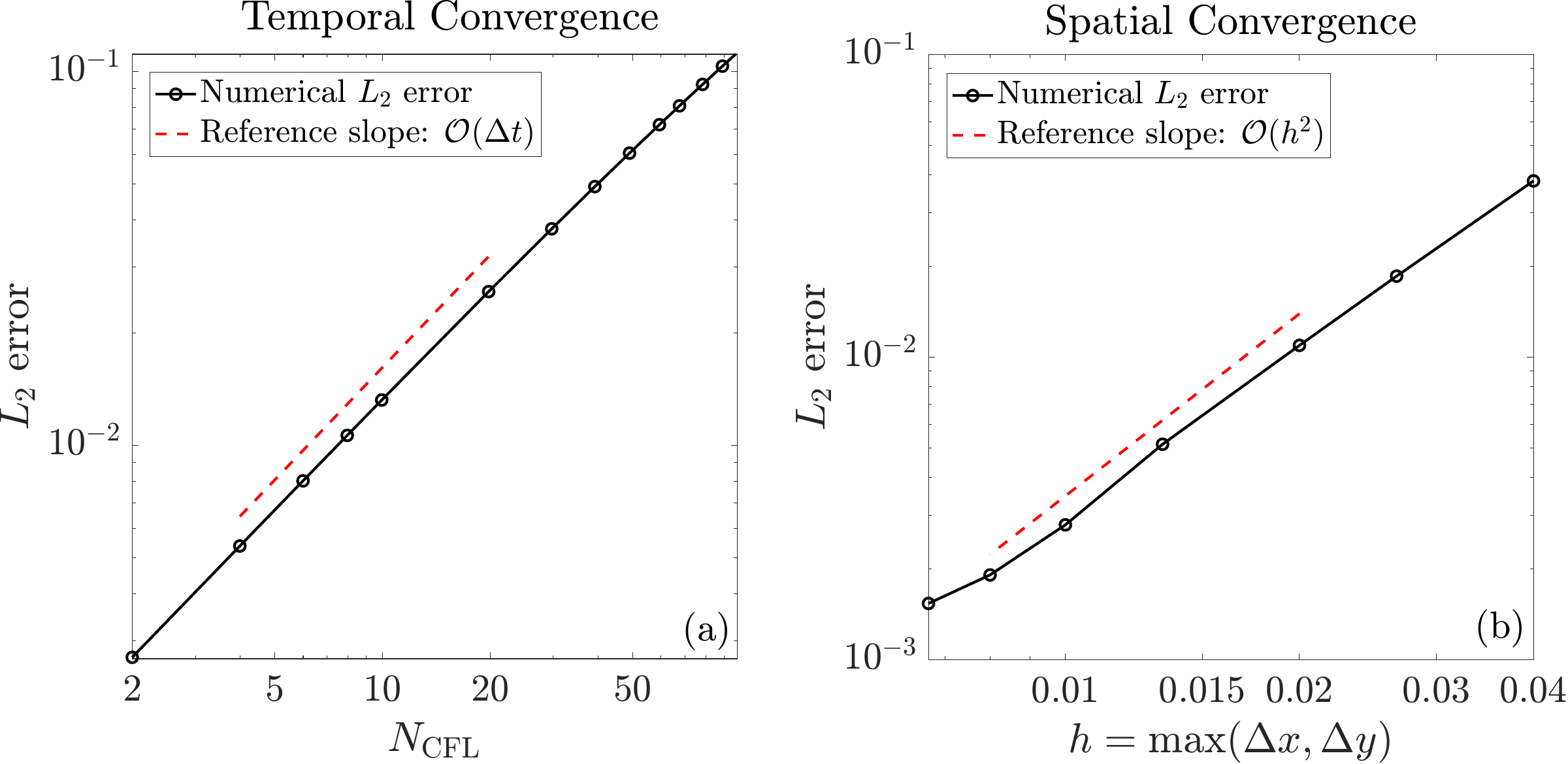}
 \caption{Log-log temporal and spatial convergence of the penalized implicit scheme for the variable-coefficient anisotropic diffusion test. (a) The $L_2$ error at the final time is plotted as a function of the time-step size $\Delta t$, together with a reference line of slope~1. (b) Spatial convergence study: the $L_2$ error at the final time is plotted as a function of the mesh size, together with the corresponding reference line of slope~2.}
    \label{fig:temp_conv}
\end{figure}

Figure~\ref{fig:temp_conv} (a) shows a temporal convergence study, depicting the $L_2$ error norm as a function of the numbers of CFL, $N_{CFL}$, which ranges from $2$ to $100$. The final time is taken as $T_f=0.1$ and the reference solution is computed using a forward Euler explicit time stepping with time step $\Delta t = 0.1\,\Delta t_{\mathrm{CFL}}$. The solution exhibits good first-order convergence; however, for large $N_{CFL}$ numbers we observe a slight degradation of the formal order of convergence, which motivates exploration of higher-order time-stepping schemes in future work. {Figure~\ref{fig:temp_conv}(b) shows a spatial convergence study, where the $L_2$ error norm is plotted as a function of the mesh spacing $h=\max(\Delta x,\Delta y)$. The reference solution is computed on a fine mesh with $N_x=N_y=600$, using a corresponding time step $\Delta t=\Delta t_{\mathrm{CFL}}=1.4\times10^{-6}$. To eliminate temporal error and isolate the spatial discretization error, the time step is held fixed in this experiment. The solution is evolved to a final time of $T_f=0.01$. The results indicate good second-order spatial convergence, with observed orders in the range $1.9$-$2.0$. The slight deterioration from second order is likely due to the SMART algorithm which can revert to local first-order spatial accuracy to preserve monotonicity \cite{chacon2025robust}.}

\begin{figure}
    \centering
    \includegraphics[width=0.8\linewidth]{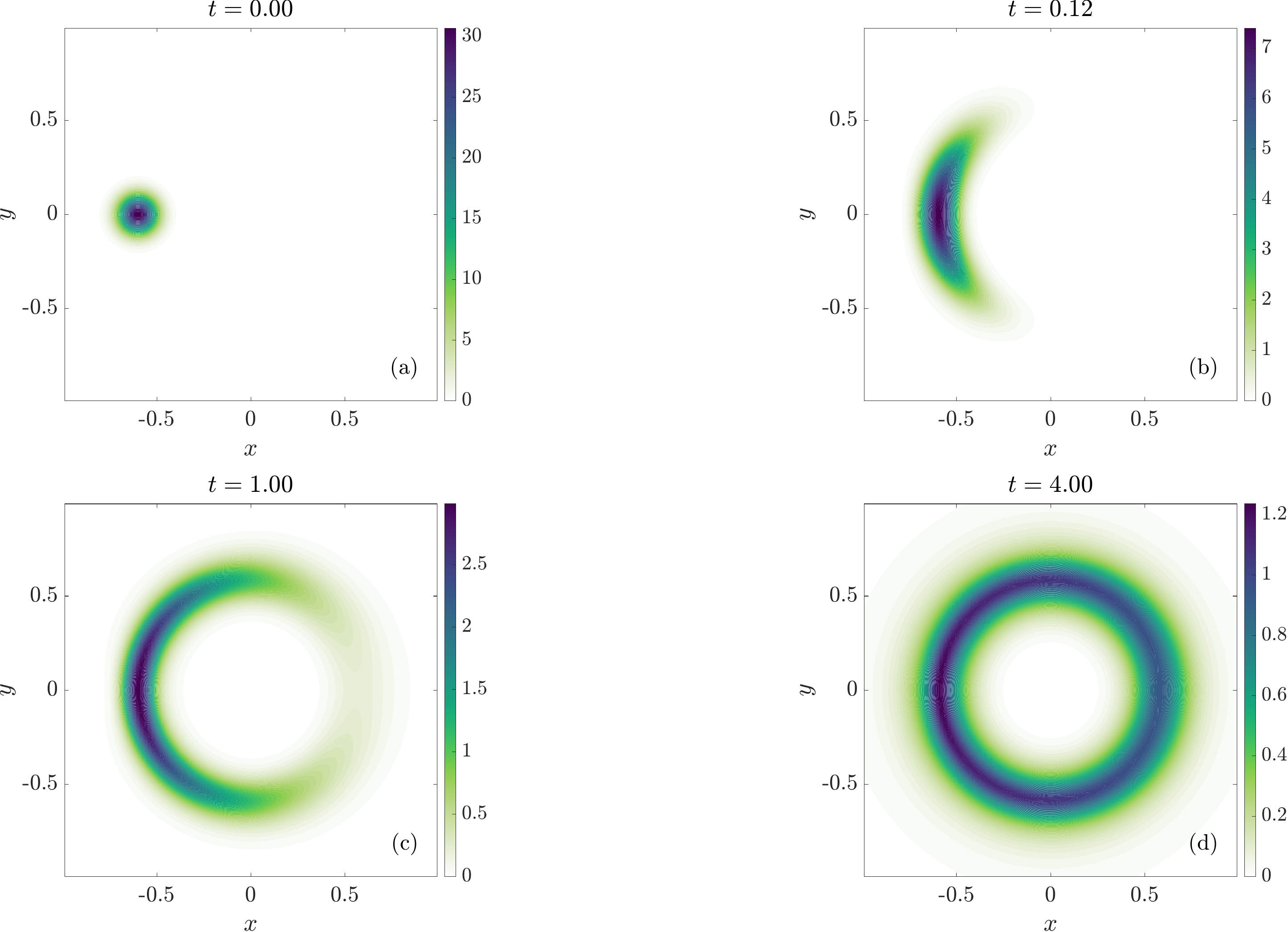}
    \caption{Evolution of the solution for the anisotropic ring test at four representative times.}
    \label{fig:ring_sim}
\end{figure}

\begin{figure}
    \centering
    \includegraphics[width=\linewidth]{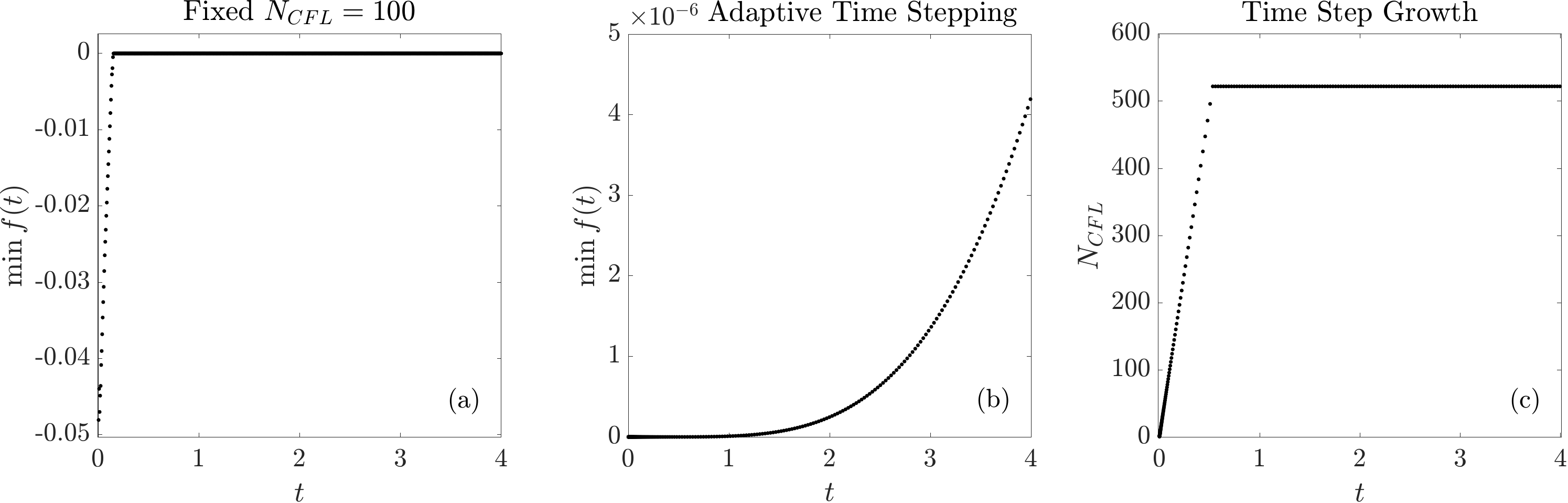}
 \caption{(a) Minimum value of the solution in time for a fixed time-step ratio with $N_{CFL}=100$. (b) Minimum value of the solution in time for the adaptive time-stepping strategy. (c) Evolution of the adaptive time-step ratio $\Delta t_n / \Delta t_{\mathrm{CFL}}$ as a function of time.}
    \label{fig:ring_negativity}
\end{figure}

Figure \ref{fig:ring_sim} shows the evolution of the initial Gaussian as function of time, as it is transported into a ring-shaped structure.
Finally, Figure~\ref{fig:ring_negativity}(a) displays the evolution of the adaptive timestep size as a function of the time-step index, showing a rapid growth as the simulation proceeds. Panels~\ref{fig:ring_negativity}(b)-(c) compares the minimum of the solution in time for a fixed time step with $N_{CFL}=100$ and for the adaptive strategy, respectively. In the fixed time-step case the solution develops pronounced negativity at early times, whereas the adaptive simulation remains positive as the time step is increased gradually up to $N_{CFL=500}$.

\subsection{Linearized Rosebluth-Fokker-Planck relaxation tests}
Here, we consider the diffusion tensor and advection coefficient form a static background species $\beta$ in its own equilibrium state \cite{liu2025optimization}. Using the equilibrium relation \eqref{eq:tensor_in_equilbirium}, the linearized RFP equation for the distribution $f_\alpha$ of a particles $\alpha$ colliding in a background medium can be written as
\begin{equation}
\label{eq:rfp_cross}
\frac{\partial f_{\alpha}}{\partial t}
=
\Gamma_{\alpha\beta}\,
\nabla_{\mathbf v}\!\cdot\!\left[
\mathcal D^{M}_{\beta}\cdot\left(
\nabla_{\mathbf v} f_{\alpha}
-\frac{m_{\alpha}}{T_{\beta}}\left(\mathbf u_{\beta}-\mathbf v\right)f_{\alpha}
\right)\right],
\end{equation}
where $\Gamma_{\alpha\beta}
=\frac{2\pi\,e_{\alpha}^{2}e_{\beta}^{2}\,\lambda\,n_{\beta}}{m_{\alpha}^{2}}$ is the collision frequency, and $\mathcal D^M_\beta=\nabla_{\mathbf v}\nabla_{\mathbf v} G^M_\beta$ is the diffusion tensor defined through the Maxwellian Rosenbluth potentials via $ \nabla_\mathbf{v}^2 G^M_\beta=H^M_\beta$ and $ \nabla_\mathbf{v}^2 H^M_\beta=f^M_\alpha$. Here $f^M_\alpha$ denotes the equilibrium Maxwellian, parameterized by the mass of the test particles $\alpha$, $m_{\alpha}$, the equilibrium drift velocity of the background species $\beta$, $\mathbf u_{\beta}$, and the equilibrium temperature of the background species $\beta$, $T_{\beta}$. Closed-form expressions for $\mathcal D^M_\beta$ can be derived directly from these potentials; for conciseness, we omit the details and refer interested readers to \cite{liu2025optimization}. {The penalization parameter is chosen as one half of the largest eigenvalue of $\mathcal{D}^M$ across the mesh.}

A few words regarding the discretization startegy of \eqref{eq:rfp_cross} are in order. Expanding the linearized operator componentwise yields isotropic contributions (corresponding to $D_{\perp\perp}$ and $D_{\parallel\parallel}$), which are of advection--diffusion type and which we discretize using the classical Chang-Cooper scheme, as detailed in Section~\ref{MP_discretization_RFP}. The anisotropic contributions, corresponding to $D_{\perp\parallel}$, are discretized using the advectionalization technique described in Section~\ref{MP_discretization_Heat}. For the temporal discretization, we use the penalization strategy developed in this manuscript. 

The linearized RFP equation serves twofold: it enables a direct comparison of the proposed numerical scheme against theoretical relaxation rates, and it highlights key aspects of the penalization framework, including equilibrium-preserving penalization and the role of variable penalization and adaptive time-stepping in maintaining positivity.

\subsubsection{Beam relaxation test}
In this section, we reproduce the test of \cite{nrl_2023,taitano2021conservative,hinton1983collisional} to assess the accuracy of the penalization scheme by simulating the linearized Rosenbluth-Fokker-Planck equation for deuterium (D) test particles with mass \(m_{D}=2\) colliding with a Maxwellian aluminum (Al) background with mass \(m_{Al}=27\). 

We compare the numerical results against the theoretical slowing-down, perpendicular (transverse) diffusion, and parallel diffusion rates, denoted by \(\nu_s^{D/Al}\), \(\nu_\perp^{D/Al}\), and \(\nu_\parallel^{D/Al}\), respectively, which are given analytically by (assuming that the colliding D-species follows a Dirac-delta distribution):
\begin{align}
\label{nus_theoretical}
\nu_s^{D/Al} &= \Bigl(1+\frac{m_{D}}{m_{Al}}\Bigr)\,\psi\!\bigl(x^{D/Al}\bigr)\,\nu_0^{D/Al},\\[0.2em]
\label{nuperp_theoretical}\nu_\perp^{D/Al} &= 2\Biggl[\Bigl(1-\frac{1}{2x^{D/Al}}\Bigr)\psi\!\bigl(x^{D/Al}\bigr)
      + \psi'\!\bigl(x^{D/Al}\bigr)\Biggr]\nu_0^{D/Al},\\[0.2em]
\label{nupar_theoretical}\nu_\parallel^{D/Al} &= \Biggl[\frac{\psi\!\bigl(x^{D/Al}\bigr)}{x^{D/Al}}\Biggr]\nu_0^{D/Al},
\end{align}
where
\begin{align}
\nu_0^{D/Al} &= \frac{4\pi\,e_{D}^2 e_{Al}^2\,\lambda\,n_{Al}}{m_{D}^{2}\,u_{\parallel,D}^{3}},
\qquad
x^{D/Al} = \frac{u_{\parallel,D}^{2}}{v_{th,Al}^{2}},\\[0.2em]
\psi(x) &= \frac{2}{\sqrt{\pi}}\int_0^x t^{1/2}e^{-t}\,dt,
\qquad
\psi'(x) = \frac{d\psi}{dx}.
\end{align}
The charges are set to \(e_{Al}=13\) and \(e_{D}=1\), the Coulomb logarithm is taken as \(\lambda=1\), and the background density is \(n_{Al}=1\).

{In our simulations, the temperatures are initialized as \(T_{Al}=m_{Al}/2=13.5\) and  \(T_{D}=0.01\) (cold enough for the Dirac-delta assumption in the theory to hold), yielding thermal velocities \(v_{th,Al}=1\) and \(v_{th,D}=0.1\) for the background aluminum and deuterium species, respectively. The initial parallel drifts are set to \(u_{\parallel,D}=30\,v_{th,Al}\) and \(u_{\parallel,Al}=0\).} To accurately resolve the initial deuterium beam, we center the computational domain around the initial drifting Maxwellian,
\[
f_0(\mathbf{v})
= \frac{1}{\pi^{3/2}\,v_{th,D}^{3}}
\exp\!\Biggl(
-\frac{(v_\parallel-u_{\parallel,D})^{2}+v_\perp^{2}}{v_{th,D}^{2}}
\Biggr),
\]
using
\[
v_\parallel \in \bigl[u_{\parallel,D}-0.5,\;u_{\parallel,D}+0.5\bigr],
\qquad
v_\perp \in \bigl[0,\;0.5\bigr],
\]
with a grid of \(N_\parallel = 128\) and \(N_\perp = 64\).
Finally, to match the theoretical time scale, we set:
\[
\tau^{0}=\frac{{v_{th,M}}^{3/2}}{2\,n_{Al}},
\qquad \text{with} \qquad 
v_{th,M}=\sqrt{\frac{2T_{Al}}{m_{D}}}.
\]
Similarly to \cite{taitano2021conservative}, the numerical slowing-down, perpendicular, and parallel diffusion rates are time-averaged to increase accuracy. They are defined by:
\begin{align}
\label{nus_numerical}
\langle \nu_s^{D/Al}\rangle_{\tau}
&= \frac{1}{N_t\,\Delta t}\sum_{p=1}^{N_t}
\left|\frac{u_{\parallel,D}^{(p)}-u_{\parallel,D}^{(p-1)}}{u_{\parallel,D}^{(p)}}\right|,
\qquad
u_{\parallel,D}^{(p)}=\frac{\langle v_\parallel, f_{D}^{(p)}\rangle_v}{\langle 1, f_{D}^{(p)}\rangle_v},
\\[0.4em]
\label{nuperp_numerical}\langle \nu_{\perp}^{D/Al}\rangle_{\tau}
&= \frac{1}{N_t\,\Delta t}\sum_{p=1}^{N_t}
\left|\frac{T_{\perp,D}^{(p)}-T_{\perp,D}^{(p-1)}}{m_{D}\,\bigl(u_{\parallel,D}^{(p)}\bigr)^2}\right|,
\qquad
T_{\perp,D}^{(p)}=\frac{m_{D}\,\langle v_{\perp}^{2}, f_{D}^{(p)}\rangle_v}{\langle 1, f_{D}^{(p)}\rangle_v}.
\\[0.4em]
\label{nuparallel_numerical}\langle \nu_{\parallel}^{D/Al}\rangle_{\tau}
&= \frac{1}{N_t\,\Delta t}\sum_{p=1}^{N_t}
\left|\frac{T_{\parallel,D}^{(p)}-T_{\parallel,D}^{(p-1)}}{m_{D}\,\bigl(u_{\parallel,D}^{(p)}\bigr)^2}\right|,
\qquad
T_{\parallel,D}^{(p)}=\frac{m_{D}\,\langle (v_\parallel-u_{\parallel,D}^{(p)})^{2}, f_{D}^{(p)}\rangle_v}{\langle 1, f_{D}^{(p)}\rangle_v}.
\end{align}
\begin{figure}
    \centering
    \includegraphics[width=0.7\linewidth]{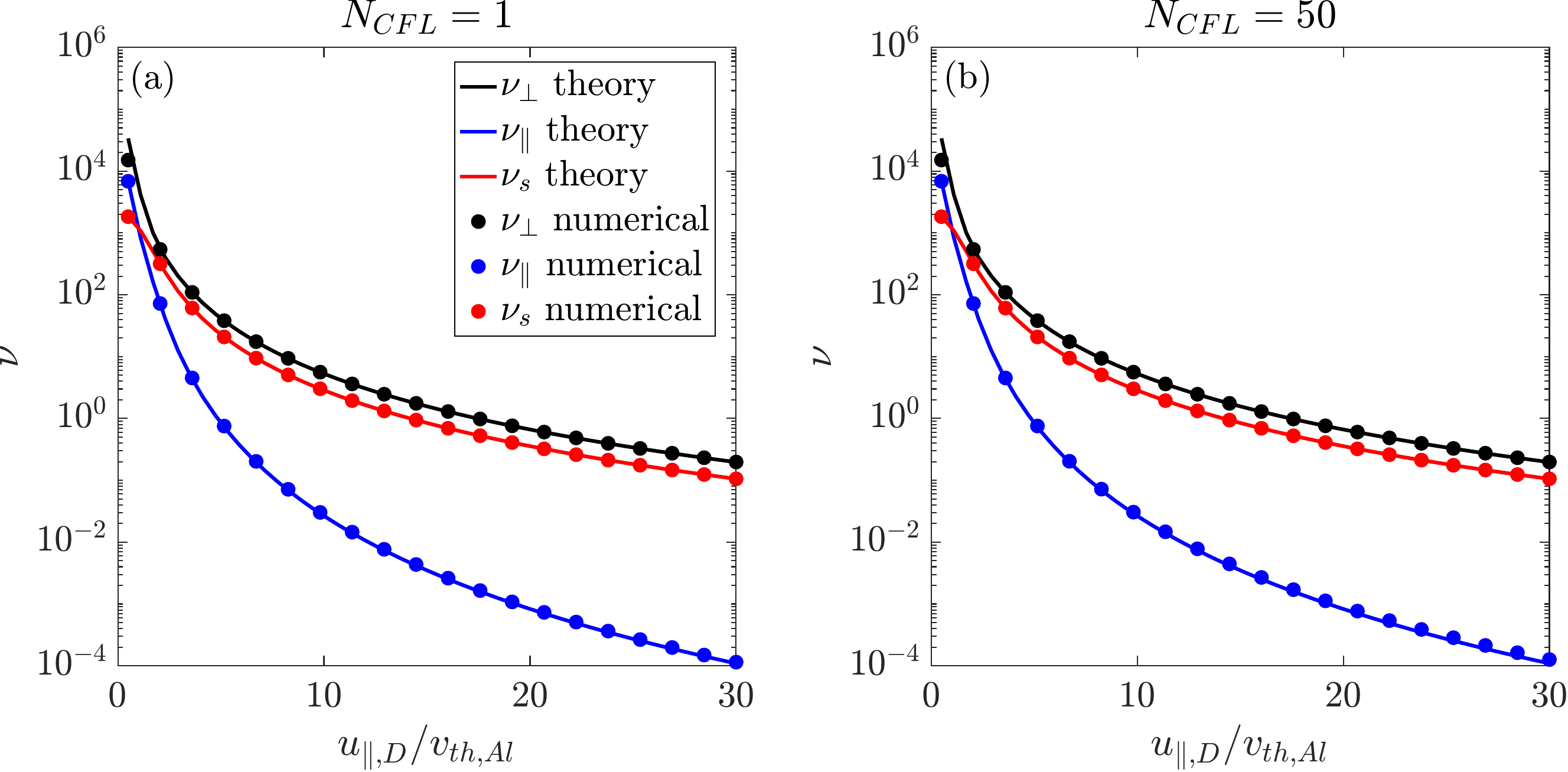}
    \caption{Beam relaxation test: Comparison between numerical (filled circles) and theoretical (solid lines) relaxation rates for a deuterium beam, including the slowing--down, perpendicular, and parallel diffusion rates. Panel (a) corresponds to a time step with $N_{\mathrm{CFL}} = 1$, while panel (b) corresponds to $N_{\mathrm{CFL}} = 50$.}
    \label{fig:Beam_Relax}
\end{figure}

{Figure~\ref{fig:Beam_Relax} shows the comparison between the theoretical rates given by \eqref{nus_theoretical}-\eqref{nupar_theoretical} and the numerical rates obtained from \eqref{nus_numerical}-\eqref{nuparallel_numerical}. Panel (a) depicts results from the penalization algorithm with an explicit time step corresponding to $N_{\mathrm{CFL}}=1$. Excellent agreement is observed across a range of drift velocities $u_\parallel$. Panel (b) depicts results from the same experiment but with a larger time step corresponding to $N_{\mathrm{CFL}}=50$, and also demonstrates very good agreement. Larger timesteps result in significant deviation of the D distribution function from the Dirac-delta in a single timestep, degrading agreement with theory.} 

\subsubsection{Cross-species pitch-angle scattering test}

In this test, we simulate the linearized RFP model ~\cite{ liu2025optimization} on the $(v_\perp,v_\parallel)$ domain $[0,6]\times[-6,6]$ with grid size $N_{v_\perp}=N_{v_\parallel}=100$. We consider collisions of an electron species with a background ion species. The initial condition for the electrons is parametrized with a normalized mass $m_e=1$, and initial temperature $T_{e,0}=0.25$ and a drift velocity $u_e=4$. The diffusion tensor is computed from a background ions species with mass $m_i=100$, drift velocity $u_i=0$, and temperature $T_b=1$. {We normalize the collision frequency by setting $\Gamma_{ei}$=1.}

\begin{figure}
    \centering
    \includegraphics[width=\linewidth]{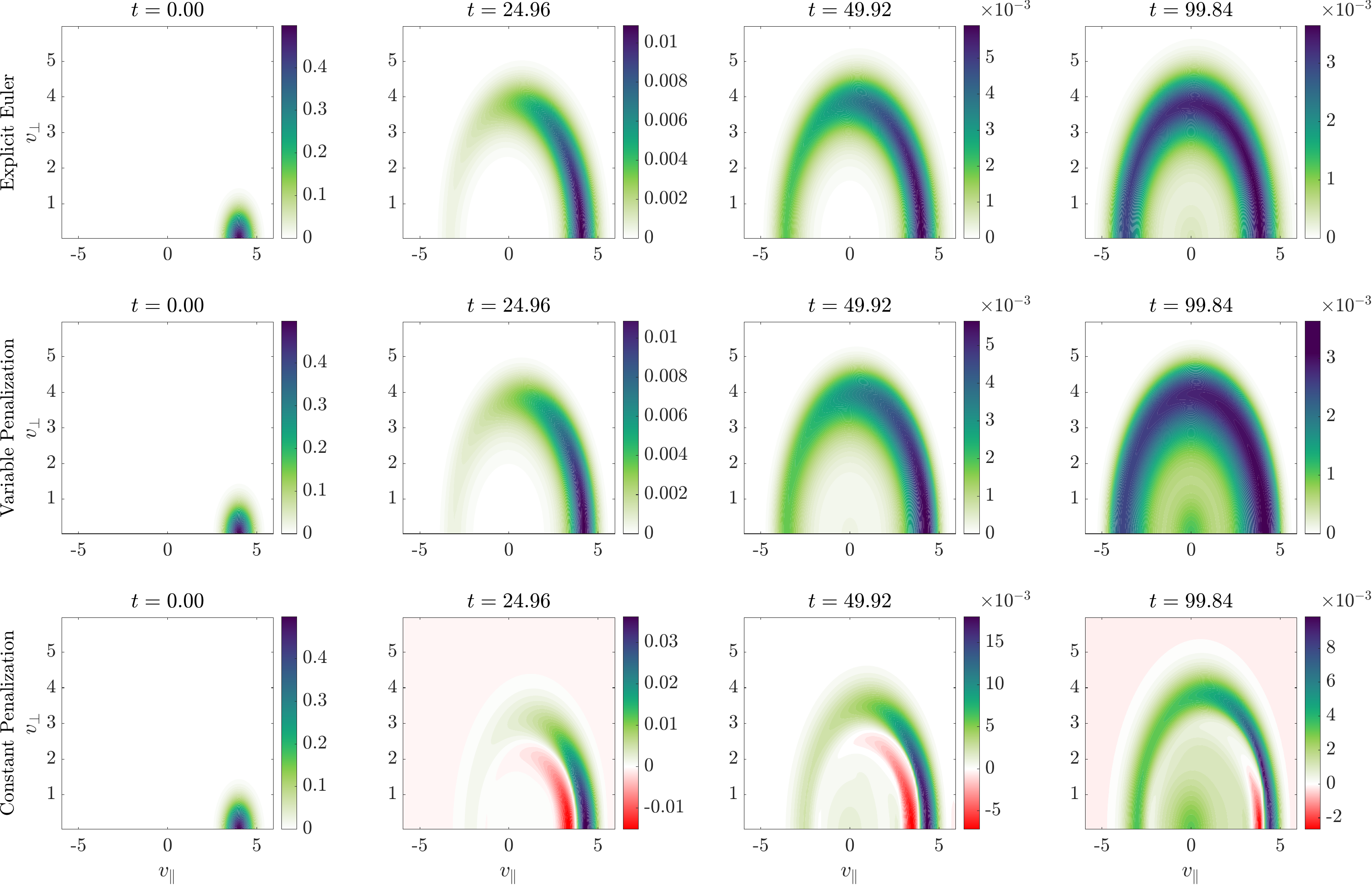}
  \caption{
    Accuracy impact of constant vs. variable penalization in pitch-angle scattering test. Rows correspond to different time-integration schemes: (top) explicit Euler, (middle) fixed-$\Delta t$ implicit scheme with variable penalization, and (bottom) fixed-$\Delta t$ implicit scheme with constant penalization. Columns show solution snapshots in the $(v_\perp,v_\parallel)$ plane at representative times $t = 0,\,24.96,\,49.92,\,99.84$, illustrating the importance of variable penalization for an accurate evolution of the distribution. 
    }
    \label{fig:Pitch_angle}
\end{figure}

We begin by looking at the accuracy impact of constant vs. variable penalization of the RFP equation with a fixed timestep. Figure~\ref{fig:Pitch_angle} shows three sets of simulations: explicit Euler time stepping (top row), penalized implicit time stepping with variable coefficients (middle row), and penalized implicit time stepping with constant penalization (bottom row). The implicit penalization schemes employ a uniform time step with $N_{CFL}=500$, with $\Delta t_{CFL} = 3.53\times 10^{-4}$. Qualitatively, the variable-coefficient penalization compares much better to the explicit scheme than the constant penalization result.

\begin{figure}
    \centering
\includegraphics[width=\linewidth]{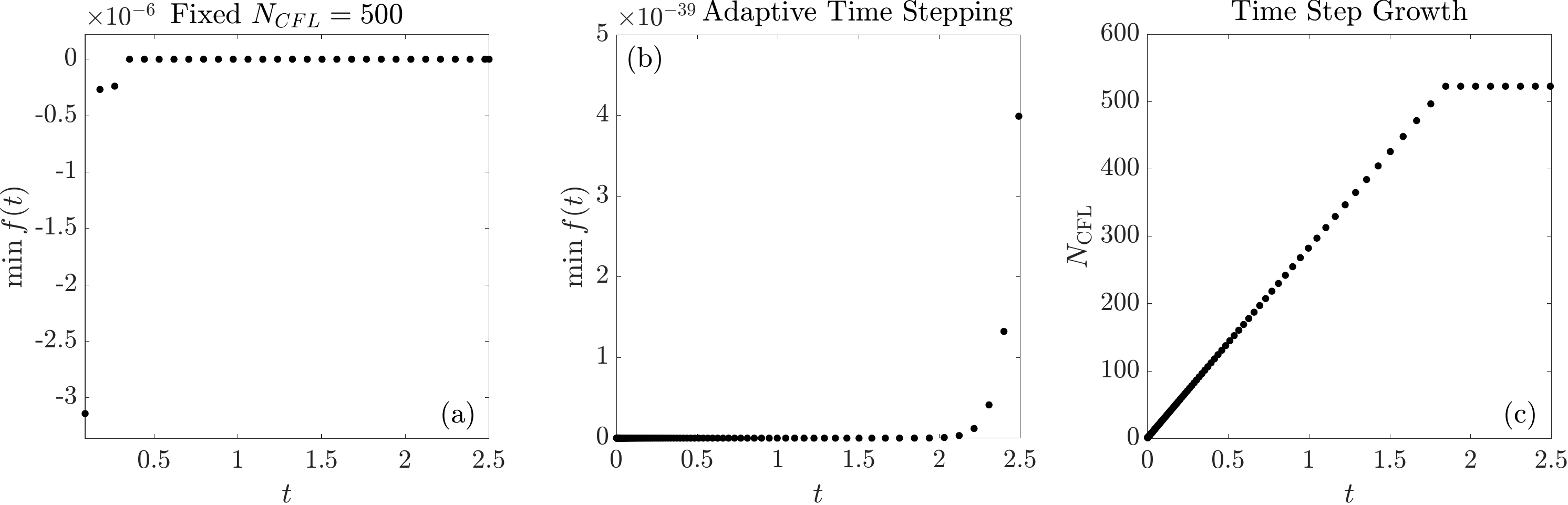}
     \caption{
    (a) Minimum value of the solution in time for a fixed time-step ratio $\Delta t = 500\,\Delta t_{\mathrm{CFL}}$. (b) Minimum value of the solution in time for the adaptive time-stepping strategy. (c) Evolution of the adaptive time-step ratio $\Delta t_n / \Delta t_{\mathrm{CFL}}$ as a function of time.}
    \label{fig:Pitch_angle_pos}
\end{figure}

Figure~\ref{fig:Pitch_angle_pos} illustrates the positivity-preserving properties of the proposed variable penalization scheme with and without variable timestepping. Panels~\ref{fig:Pitch_angle_pos}(a) and~\ref{fig:Pitch_angle_pos}(b) compare constant time stepping with a variable time-stepping strategy, respectively. We clearly observe that the constant time-stepping approach exhibits negativity at the beginning of the simulation on the order of $10^{-6}$, whereas adaptive time-stepping strategy recovers and maintains positivity. Panel~\ref{fig:Pitch_angle_pos}(c) shows the growth of the adaptive $N_{CFL}$ as a function of the time-step iteration index (corresponding to the simulation in Panel~\ref{fig:Pitch_angle_pos}(b)), displaying the rapid increase of the time-step using the proposed prescription. The time step grows exponentially until it reaches the target value
$N_{\mathrm{CFL},\max} \ge 500$, at which point it is fixed to match the
uniform time-stepping run with
$\Delta t = N_{\mathrm{CFL},\max}\,\Delta t_{\mathrm{CFL}}$.

\begin{figure}
    \centering
\includegraphics[width=0.6\linewidth]{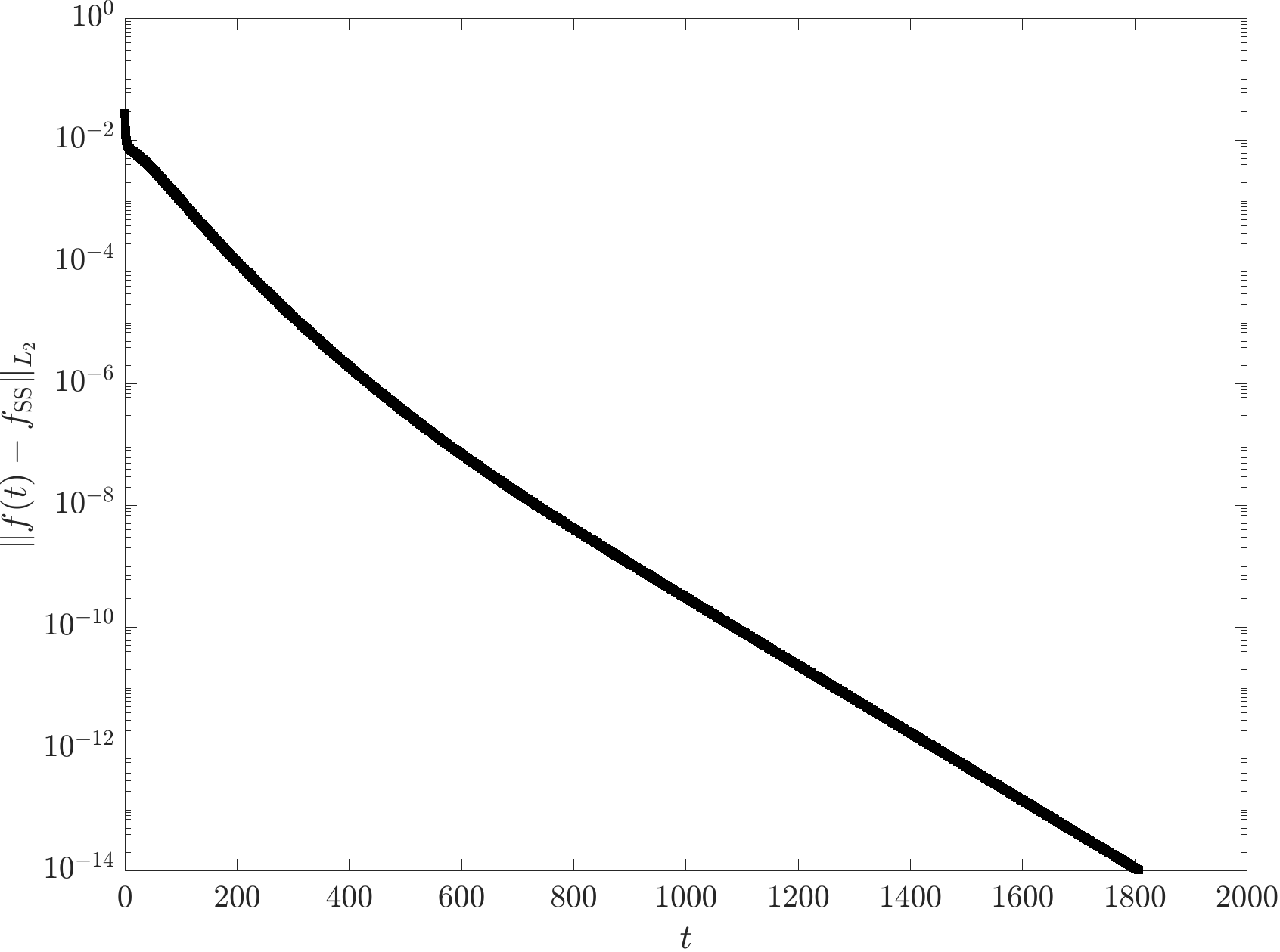}
 \caption{Semi-logarithmic plot of the convergence to the steady-state solution, denoted by $f_{\mathrm{SS}}$ in the figure, for the pitch-angle scattering test. The $L_2$-norm of the error between the variable-penalization adaptive-time solution and the steady-state solution is computed at each time step.}
    \label{fig:Pitch_angle_conv}
\end{figure}

Figure~\ref{fig:Pitch_angle_conv} illustrates the convergence of the solution to the steady-state Maxwellian up to machine precision. For this test, we retain all parameters from the previous setup, except for the background mass, which is set to $m_b = 5$ to obtain a more tractable convergence. This experiment confirms the steady-state preservation property of the proposed algorithm.

\subsection{Nonlinear Rosenbluth-Fokker-Planck relaxation tests}

\subsubsection{Temperature isotropization test}
To verify the nonlinear penalized RFP scheme, we consider a classical temperature isotropization
test.  We consider different initial perpendicular
and parallel temperatures, $T_\perp(0)=0.5$ and $T_\parallel(0)=0.1$.  Collisions drive these temperatures toward an equilibrium temperature
according to:
\begin{equation}
    \frac{dT_\perp}{dt}
    =
    -\frac{1}{2}\,\frac{dT_\parallel}{dt}
    =
    -\nu_T\bigl(T_\perp - T_\parallel\bigr),
\end{equation}
where the isotropization frequency $\nu_T$ is given by

\begin{equation}
\nu_T
=
\frac{2\sqrt{\pi}\, e^{4}\,n\,\lambda}
     {m^{1/2}\, T_{\parallel}^{3/2}}
\,A^{-2}
\left[
  -3 + (A+3)\,\frac{\tan^{-1}\!\bigl(A^{1/2}\bigr)}{A^{1/2}}
\right],
\qquad
A \equiv \frac{T_\perp}{T_\parallel}-1 .
\end{equation}
with the plasma parameters (charge $e$ , density $n$, masses $m$, and Coulomb
logarithm $\lambda$), see \cite{nrl_2023}.

To compare against the semi-analytical theoretical predictions, we consider the domain $[0,5]\times[-5,5]$ in the $v_\perp$ and $v_\parallel$ directions, respectively, with grid sizes $N_{v_\perp}$ and $N_{v_\parallel}$ satisfying $N_{v_\perp} = N_{v_\parallel}/2 = 64$. We consider the initial condition:

\[
    f_0(v_\perp, v_\parallel)
    =
    \frac{1}{2^{3/2}\,\pi^{3/2}\,v_{th,\perp}^2\,v_{th,\parallel}}
    \exp\!\left(
        -\frac{v_\perp^2}{2\,v_{th,\perp}^2}
        -\frac{(v_\parallel - 0.05)^2}{2\,v_{th,\parallel}^2}
    \right),
\]
where the thermal speeds are defined as:
\[
    v_{th,\perp} = \sqrt{\frac{T_\perp}{m}},
    \qquad
    v_{th,\parallel} = \sqrt{\frac{T_\parallel}{m}}.
\]

\begin{figure}
    \centering
    \includegraphics[width=0.6\linewidth]{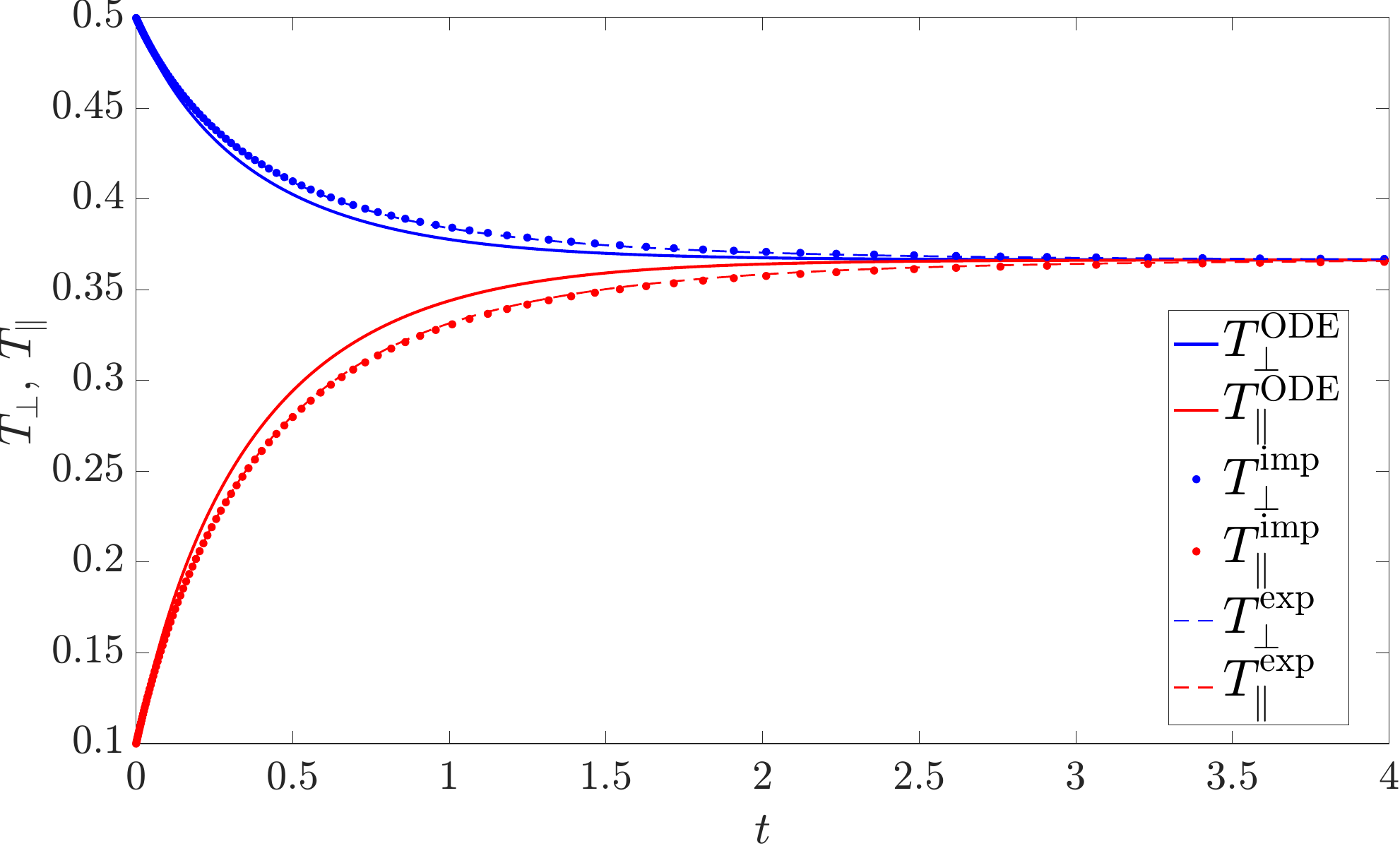}
    \caption{Isotropization test of perpendicular and parallel temperatures, $T_\perp$ and $T_\parallel$, respectively. Three simulations are shown: (i) semi-implicit model (solid lines), (ii) implicit penalized scheme (circles), and (iii) explicit scheme (dashed lines). 
    }
    \label{fig:Isotropization}
\end{figure}

Figure~\ref{fig:Isotropization} shows the evolution of the parallel and perpendicular temperatures as a function of time. We observe excellent agreement between the proposed penalized algorithm and a reference solution obtained with an explicit Euler scheme. Both solutions also show good agreement with the approximate semi-analytical solution presented above, integrated with MATLAB's ODE solver \texttt{ode23}.

\begin{figure}
    \centering
    \includegraphics[width=\linewidth]{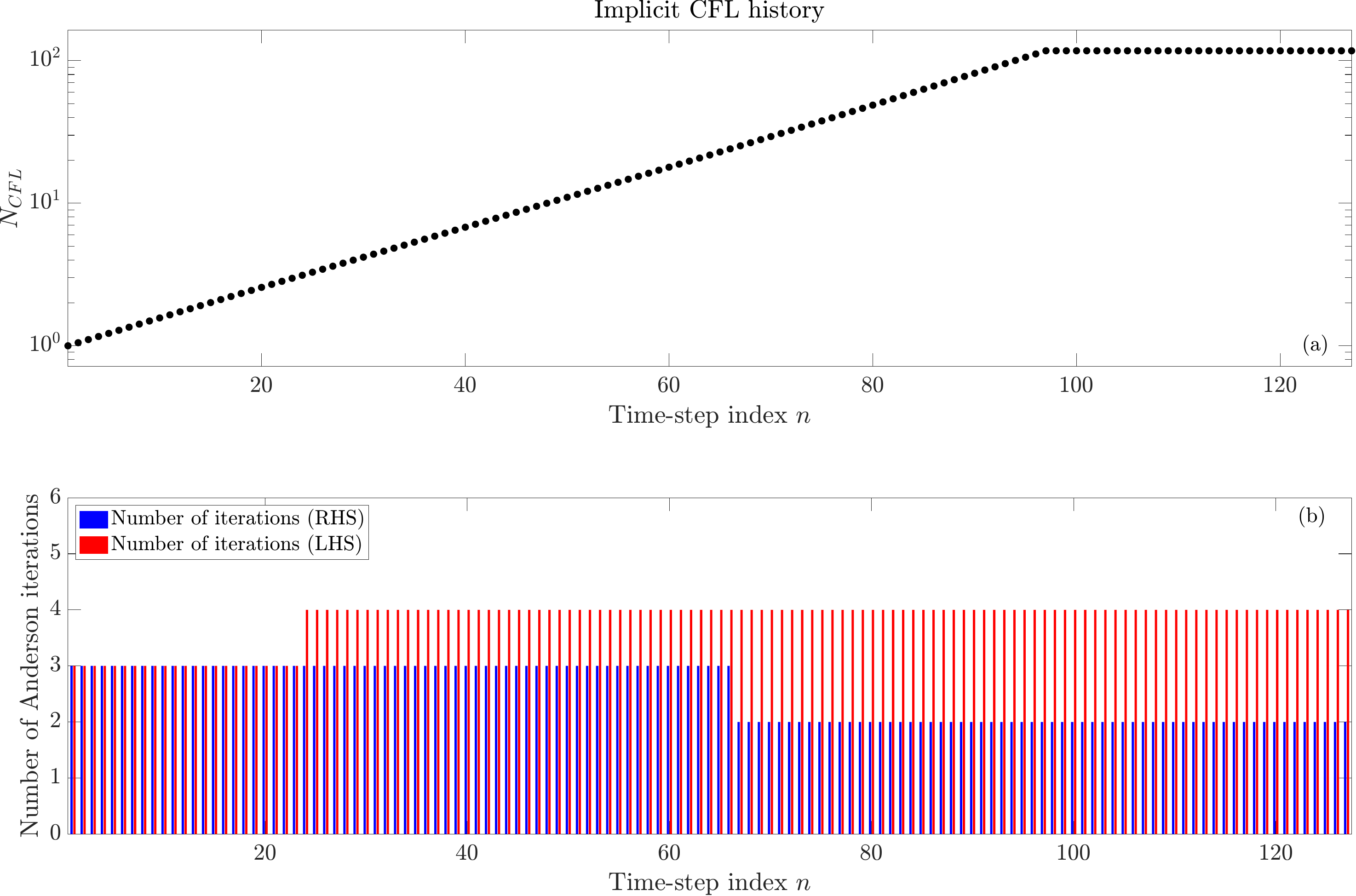}
    \caption{Isotropization test: (a) evolution of the time-step growth under time adaptivity; (b) number of Anderson-acceleration iterations for the Lagrange multiplier associated with the left-hand-side penalization operator (blue) and the right-hand-side operator (red). }
    \label{fig:CFL_fig_isotrop}
\end{figure}
Figure~\ref{fig:CFL_fig_isotrop}(a) shows the evolution of the adaptive timestep as a function of the time-step index. The time step increases according to our adaptive time-stepping strategy in an exponential fashion until it reaches the prescribed maximum CFL of $100$. Figure~\ref{fig:CFL_fig_isotrop}(b) shows the number of Anderson-acceleration nonlinear iterations required for enforcing conservation, which is essentially constant across time steps and does not exceed 4 iterations for the $(u_{\parallel,\beta}^{n+1},\lambda_\beta^{n+1})$ solve, and 3 iterations for the $(u_{\parallel,\beta}^{n},\lambda_\beta^{n})$ one.
\begin{figure}
    \centering
    \includegraphics[width=\linewidth]{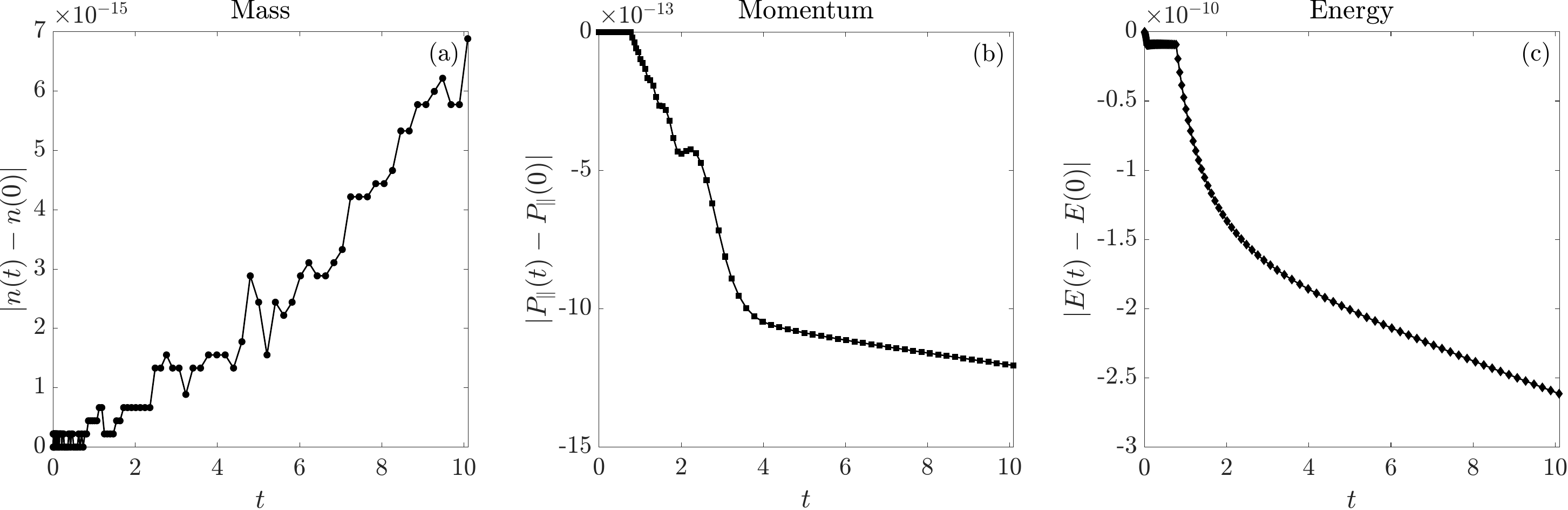}
    \caption{History of conservation error in the moments: (a) absolute density, (b) absolute momentum, and (c) absolute energy. The absolute error is computed with respect to the initial moments at time $t=0$.}
    \label{fig:cons_isotrop}
\end{figure}
Figures~\ref{fig:cons_isotrop}(a)-(c) show the conservation properties of the solution for mass, momentum, and energy. Mass is conserved to machine precision, and momentum and energy are conserved to within the absolute nonlinear tolerance prescribed~($10^{-10}$).

\subsubsection{Multi-Gaussian Relaxation test}
This test considers the simultaneous relaxation of multiple Gaussians. 
We again consider a $[0,5]\times[-5,5]$ velocity domain in $(v_\perp,v_\parallel)$, 
discretized with $N_{v_\perp}=64$ and $N_{v_\parallel}=128$, so that 
$N_{v_\perp} = N_{v_\parallel}/2$.
We start with a bi-Maxwellian initial condition of the form
\begin{equation}
    f_0(v_\perp,v_\parallel)
    =
    \sum_{k=1}^2
    \frac{n_k}{(2\pi)^{3/2} v_{th,0}^3}
    \exp\!\left(
        -\frac{v_\perp^2}{2 v_{th,0}^2}
        -\frac{(v_\parallel - u_{\parallel,k})^2}{2 v_{th,0}^2}
    \right),
\end{equation}
where the reference temperature and mass are fixed to $T_0 = 0.1$ and $m=1$, 
with the thermal speed
\begin{equation}
    v_{th,0} = \sqrt{\frac{T_0}{m}}.
\end{equation}

\begin{figure}
    \centering
    \includegraphics[width=0.8\linewidth]{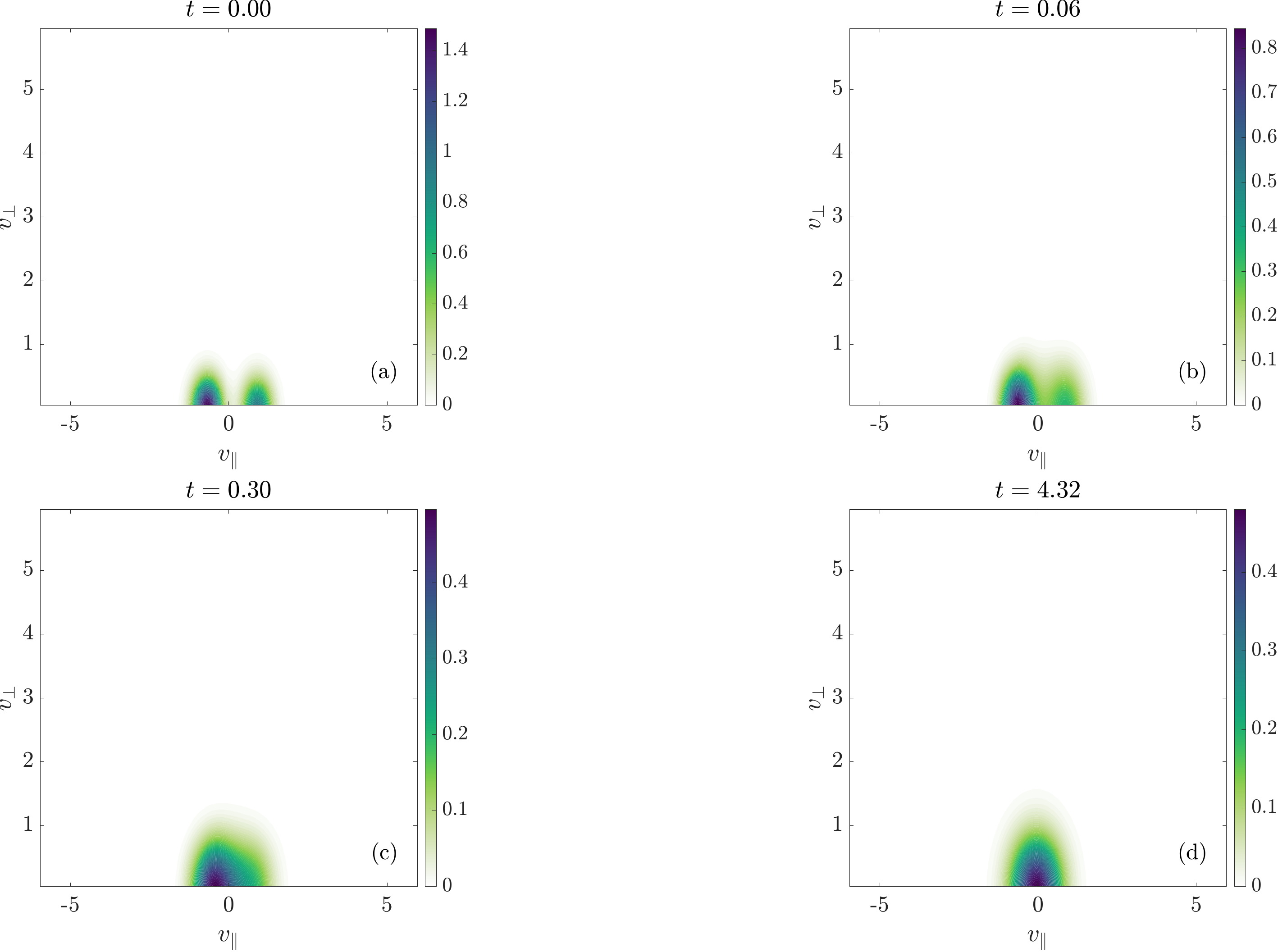}
    \caption{Evolution of the solution for the multiple Gaussian relaxation test at four representative times.}
    \label{fig:snapshots_gaussians}
\end{figure}

Figure~\ref{fig:snapshots_gaussians} shows snapshots of the solution at four representative times, 
$t = 0$, $t = 0.06$, $t = 0.30$, and $t = 4.32$. Starting from a bi-Maxwellian initial 
condition, the distribution undergoes drift and collisional diffusion, and relaxes 
toward a steady-state Maxwellian characterized by the conserved parallel drift 
$u_\parallel$ and conserved thermal velocity $v_{th,M}$.

\begin{figure}
    \centering
    \includegraphics[width=0.8\linewidth]{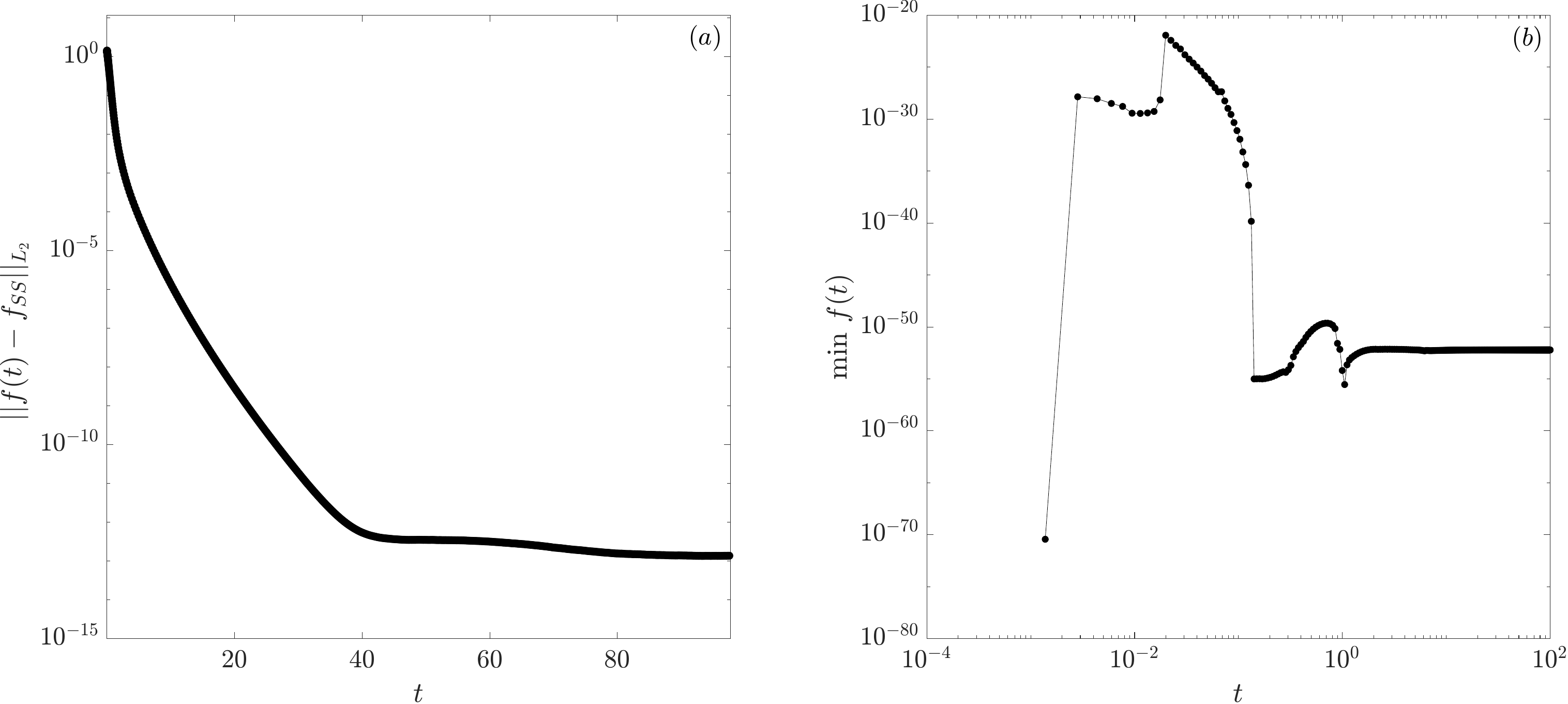}
  \caption{(a) Semi-logarithmic convergence to the steady-state Maxwellian solution. (b) Log-log plot of the minimum solution value on the grid, illustrating the positivity preserving property of the scheme.}
    \label{fig:conv_neg}
\end{figure}

Figure~\ref{fig:conv_neg}(a) shows the convergence of the numerical solution toward the analytical steady-state Maxwellian. Due to the modified Chang-Cooper discretization and the discrete conservation of mass, momentum, and energy, the normed difference between the solution and steady-state Maxwellian decays to machine precision indicating exact convergence to the analytical steady state. Figure~\ref{fig:conv_neg}(b) shows the minimum of the numerical solution in time; we observe that, under the time-adaptive strategy, the solution remains nonnegative to within machine precision throughout the simulation.

\begin{figure}
    \centering
    \includegraphics[width=\linewidth]{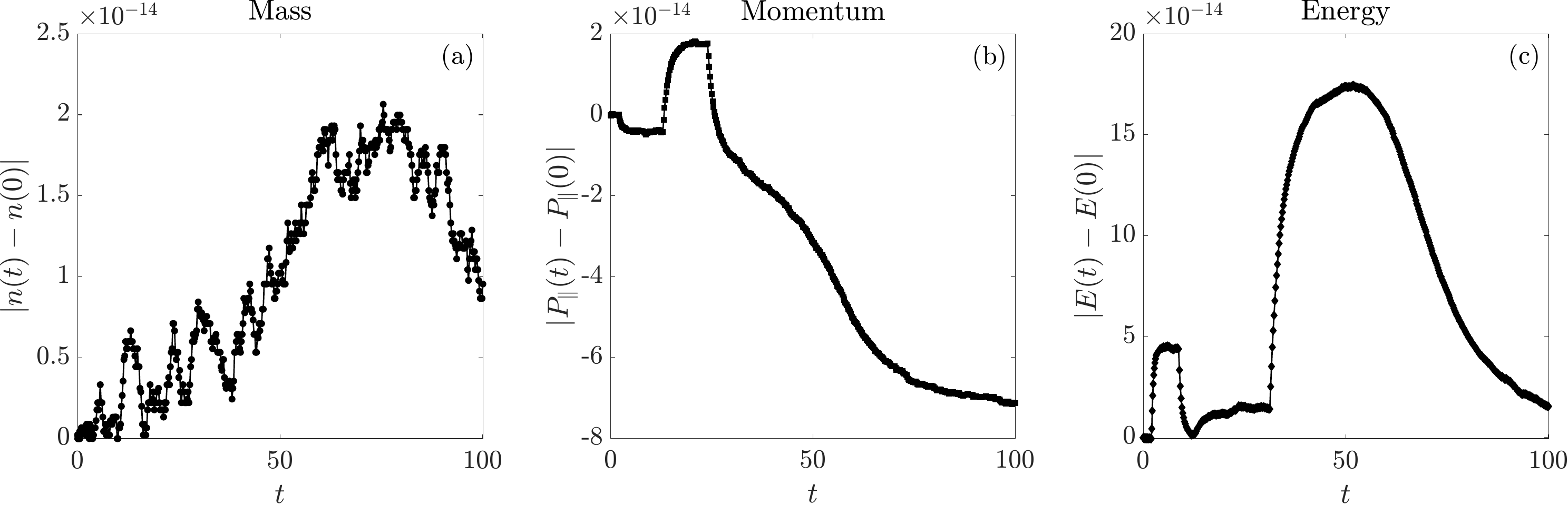}
    \caption{History of conservation error in the moments: (a) absolute density, (b) absolute momentum, and (c) absolute energy. The absolute error is computed with respect to the initial moments at time $t=0$.}
    \label{fig:2_gaussians_cons}
\end{figure}

Figure~\ref{fig:2_gaussians_cons} shows the absolute change in mass, momentum, and energy as a 
function of time. For this test, we tighten the Anderson-acceleration tolerance to 
$10^{-14}$, and we observe that the conservation errors plateau at this level. This 
confirms the excellent conservation properties of the proposed scheme.

\begin{figure}
    \centering
    \includegraphics[width=0.85\linewidth]{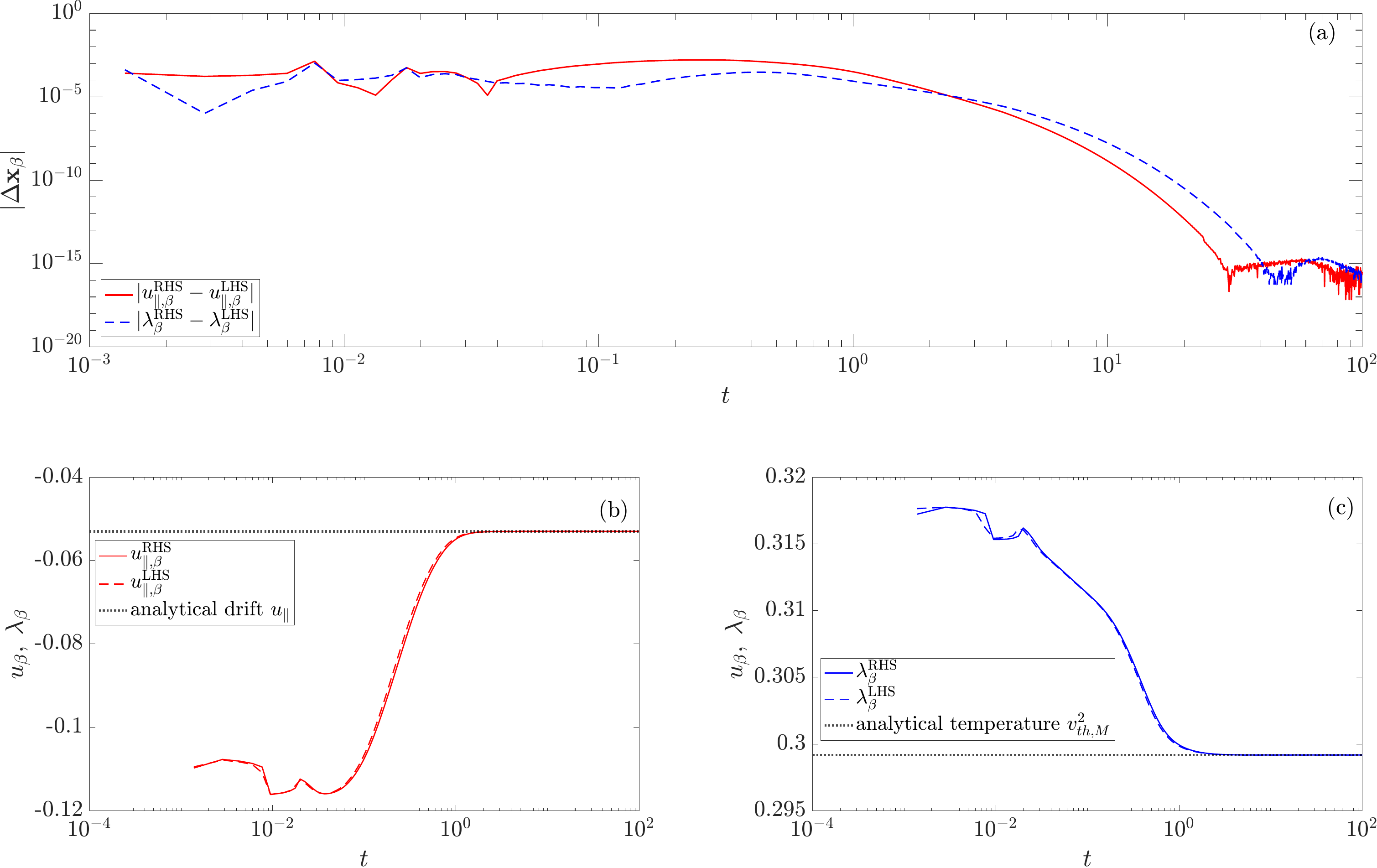}
    
\caption{Moment convergence for the conservative penalized update. 
(a) Anderson-Acceleration residuals
$\lvert u_{\parallel,\beta}^{\mathrm{RHS}} - u_{\parallel,\beta}^{\mathrm{LHS}} \rvert$ and
$\lvert \lambda_{\beta}^{\mathrm{RHS}} - \lambda_{\beta}^{\mathrm{LHS}} \rvert$
in log-log scale.
(b) Evolution of the drift velocity
$u_{\parallel,\beta}^{\mathrm{RHS}}$ and $u_{\parallel,\beta}^{\mathrm{LHS}}$ together with the
analytical equilibrium drift $u_\parallel$.
(c) Evolution of the temperature multipliers
$\lambda_{\beta}^{\mathrm{RHS}}$ and $\lambda_{\beta}^{\mathrm{LHS}}$ together
with the analytical equilibrium temperature $v_{th,M}^2$.}
    \label{fig:lagrange_moments}
\end{figure}

Figure~\ref{fig:lagrange_moments}(a)-(b) shows the difference in $(u_{\parallel,\beta},\lambda_\beta)$ between subsequent timesteps. At early times in the simulation, the difference is of order $10^{-5}$. As the simulation progresses, these moments on both 
sides converge to the same values, and their difference decreases to machine precision as the solution reaches steady state. Figure~\ref{fig:lagrange_moments}(c) shows the individual convergence of 
the Lagrange moments $u^{n}_{\parallel,\beta}$ and $u^{n+1}_{\parallel,\beta}$ toward the steady-state drift 
$u_\parallel$, confirming the steady-state preservation property of our algorithm. A 
similar trend is observed for the thermal velocities $\lambda^{n}_\beta$ and 
$\lambda^{n+1}_\beta$, which converge to the equilibrium thermal velocity $v_{th,M}^2$.

\begin{figure}
    \centering
    \includegraphics[width=0.85\linewidth]{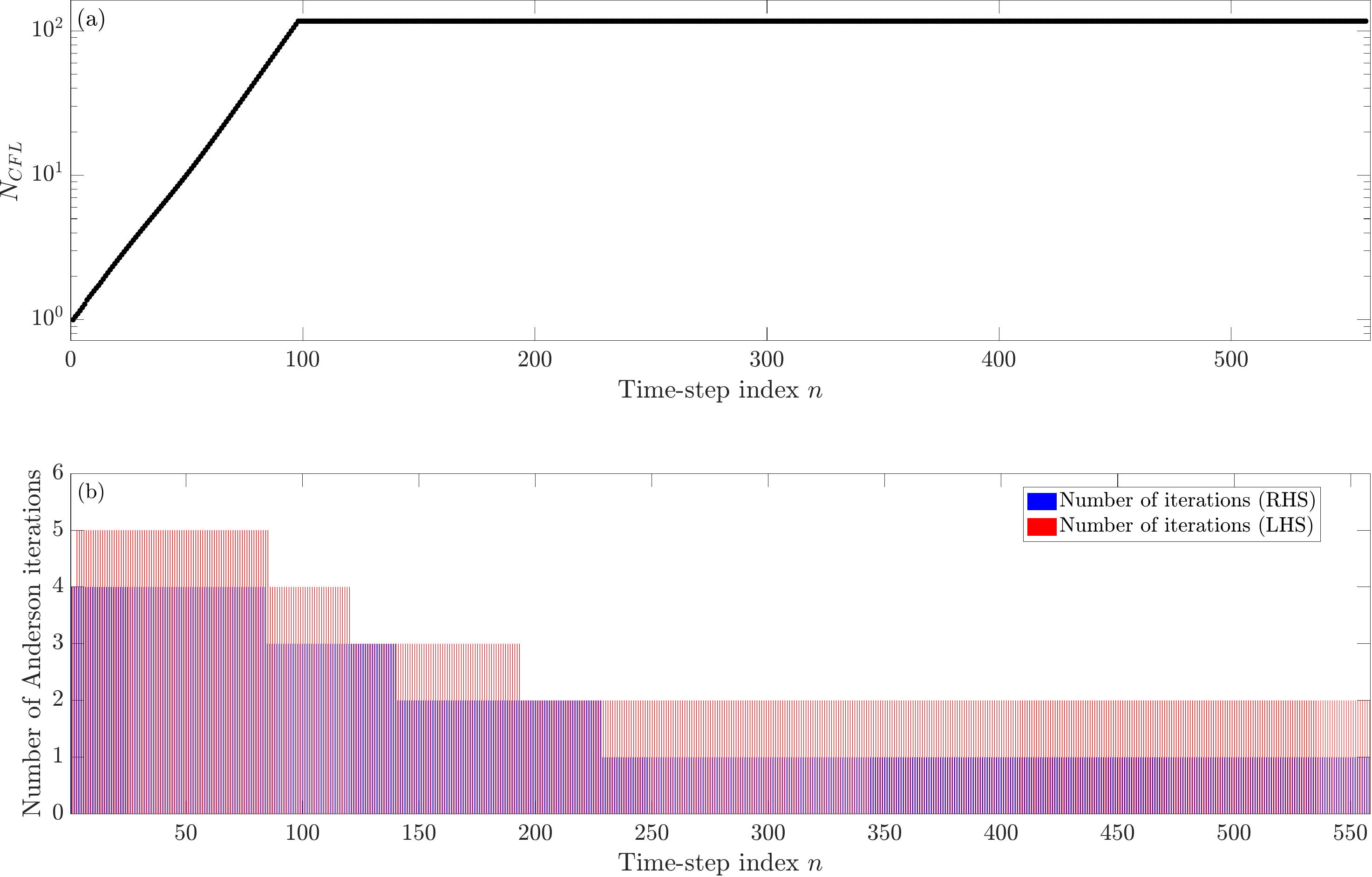}
    \caption{(a) Evolution of the adaptive time-step size; (b) number of Anderson-acceleration iterations for the left-hand-side penalization operator solve (blue) and the right-hand-side operator solve (red).}
    \label{fig:CFL_fig_isotrop_test2}
\end{figure}

Figure~\ref{fig:CFL_fig_isotrop_test2}(a) shows the history of the adaptive time-step size (on a logarithmic scale) as a function of the time-step index. The time-step growth is halted once $\Delta t$ reaches $0.2$, which corresponds to approximately $N_{CFL}=116$. 
Figure~\ref{fig:CFL_fig_isotrop_test2}(b) displays the number of Anderson-acceleration iterations required to enforce conservation of the penalization operator on the left-hand and right-hand sides. The left-hand-side iterations attain a maximum of five at early times and decrease to two as the system approaches equilibrium. Similarly, the right-hand-side iterations start around four and decay to a single iteration near equilibrium. The large admissible time steps, together with the bounded number of iterations needed to ensure conservation, attest to the efficiency of the proposed scheme.

\section{Conclusions}
\label{sec:conclusions}

We have proposed a fully discrete structure-preserving penalization scheme for the single-species RFP equation that 1) features strict conservation properties, 2) is equilibrium preserving (i.e., the analytical Maxwellian is in its null space), and 3) is able to use timesteps much larger than explicit CFLs efficiently {while robustly preserving positivity}. The penalization operator has an isotropic advection-diffusion form, and allows using timesteps much larger than explicit ones. Both the penalization and the RFP operators are strictly equilibrium preserving and feature a discrete maximum principle. In the RFP case, these properties stem from a novel generalization of the well-known Chang-Cooper discretization \cite{chang1970practical} for the RFP operator, combined with a strictly conservative discretization strategy \cite{taitano2015mass}. The discrete system is first-order accurate in time, second-order in space, and mildly nonlinear to ensure strict overall conservation. We propose a simple iterative scheme that converges very quickly and only requires linear advection-diffusion implicit solves. We have demonstrated the superior properties of the scheme with several numerical examples of increasing complexity. Future work will extend this approach to multiple species and to couple it with the Vlasov equation to study spatial transport.

\section*{Acknowledgments}
This work was performed under the auspices of the U.S.\ Department of Energy by Los Alamos National Laboratory under contract DE-AC52-06NA25396 and supported by the CHaRMNET Mathematical Multifaceted Integrated Capability Center of the Office of Applied Scientific Computing Research. JQ and JH were also supported by DoD MURI FA9550-24-1-0254. JH was additionally supported by NSF DMS-2409858.
%and DoD MURI FA9550-24-1-0254. 

\appendix

\section{Positivity of $\lambda_\beta$ from \eqref{eq:nonlinear_moments}}
\label{app:lb_moms_well_posedness}

The moment $\lambda_\beta$ is given by:
\begin{equation}
    \label{eq:nonlinear_moments_2}
    \lambda_\beta [f]
    = \frac{n_\beta E_\beta - \bigl(\boldsymbol{p}_\beta\bigr)^2}
           {n_\beta B_\beta - {\boldsymbol{A}}_\beta \cdot \boldsymbol{p}_\beta},
\end{equation}
with
\begin{equation}
    {\boldsymbol{A}}_\beta := -\int_\Omega \beta \,\nabla_\mathbf{v} f\,d\mathbf{v}, 
    \qquad
    \boldsymbol{p}_\beta := \int_\Omega \mathbf{v}\,\beta f\,d\mathbf{v},
\end{equation}
\begin{equation}
    n_\beta := \int_\Omega \beta f\,d\mathbf{v},
    \qquad
    B_\beta := -\int_\Omega \,\beta \mathbf{v}\cdot\nabla_\mathbf{v} f\,d\mathbf{v},
    \qquad
    E_\beta := \int_\Omega \beta\,v^2 f\,d\mathbf{v}.
\end{equation}
{For stability, it is required that $0 < \lambda_\beta < \infty$, as it is an effective temperature.} We show next that this is the case for arbitrary $f$.

We begin by showing that the numerator of $\lambda_\beta$ in \eqref{eq:nonlinear_moments_2} is positive for PDFs with finite support. This result follows directly from the following identity:
\begin{equation}
\label{eq:lambda_num}
0 \leq \int_\Omega d\mathbf{v} \int_\Omega d\mathbf{v}' (f\beta)(f\beta)' \, (\mathbf{v} - \mathbf{v}')^{2}
= \int_\Omega d\mathbf{v} \int_\Omega d\mathbf{v}' \, (f\beta)(f\beta)' \, (v^{2} + (v')^{2} - 2\,\mathbf{v}\cdot\mathbf{v}')
= 2\big[\, n_{\beta} E_{\beta} - \boldsymbol{p}_{\beta}\cdot\boldsymbol{p}_{\beta} \,\big],
\end{equation}
where $\beta=\beta(\mathbf{v})$ and $\beta'=\beta(\mathbf{v}')$, and similarly with $f$ and $f'$, and in the last step we have used the definitions above. Equality occurs only if $f$ is a Dirac-delta function; otherwise \eqref{eq:lambda_num} is a strict inequality for any $f$ with finite support. 

Therefore, positive and finite $\lambda_\beta$ requires the denominator in \eqref{eq:nonlinear_moments_2} to be strictly positive. To show this in general for PDFs with finite support, we begin by rewriting the denominator as:
\begin{equation}
\label{eq:lambda_den}
    n_\beta B_\beta - {\boldsymbol{A}}_\beta \cdot \boldsymbol{p}_\beta = \frac{1}{2} \int_\Omega d\mathbf{v} \int_\Omega d\mathbf{v}' \beta \, \beta' (\mathbf{v}-\mathbf{v}')\cdot (\nabla_{v'} - \nabla_\mathbf{v})f \, f'.
\end{equation}
We show next that the right hand side is positive definite for a multivariate Gaussian, and by extension to a collection to multivariate Gaussians (which can approximate continuous and positive distribution functions that vanish at infinity with arbitrary accuracy per Theorem 5 in \cite{nguyen2020approximation}). Consider first a single multivariate Gaussian:
\begin{equation}
    f_g(\mathbf v)= C \exp \left ( -\frac{1}{2} (\mathbf{v} - \mathbf{u}_g)\cdot  \Sigma^{-1}_g \cdot (\mathbf{v} - \mathbf{u}_g)\right),
\end{equation}
with $\mathbf{u}_g$ its mean, and $\Sigma_g$ its regular, symmetric-positive-definite (SPD) covariance matrix. It then follows that:
\begin{equation}
    \nabla_\mathbf{v} f_g(\mathbf{v})=- \Sigma_g^{-1}\cdot (\mathbf{v} - \mathbf{u}_g) f(\mathbf{v}).
\end{equation}
Introducing this result in \eqref{eq:lambda_den} above gives:
\begin{equation}
    \int_\Omega d\mathbf{v} \int_\Omega d\mathbf{v}' \beta \, \beta' (\mathbf{v}-\mathbf{v}')\cdot (\nabla_{v'} - \nabla_\mathbf{v})f_g \, f_g' = \int_\Omega d\mathbf{v} \int_\Omega d\mathbf{v}' \beta \, \beta' (\mathbf{v} - \mathbf{v}')\cdot  \Sigma^{-1}_g \cdot (\mathbf{v} - \mathbf{v}') f_g \, f_g' > 0,
\end{equation}
where the inequality follows from the SPD property of the covariance matrix.
The result generalizes straightforwardly for a (possibly infinite) collection of Gaussians approximating a continuous PDF $f$ with arbitrary accuracy, namely:
\begin{equation}
    f(\mathbf{v}) = \sum_g a_g f_g (\mathbf v ),
\end{equation}
with $a_g>0$. In this case, \eqref{eq:lambda_den} gives:
\begin{equation}
\label{eq:den_pos_proof}
    \int_\Omega d\mathbf{v} \int_\Omega d\mathbf{v}' \beta \, \beta' (\mathbf{v}-\mathbf{v}')\cdot (\nabla_{v'} - \nabla_\mathbf{v})f \, f' = \sum_{g,g'} a_g a_{g'} \int_\Omega d\mathbf{v} \int_\Omega d\mathbf{v}' \beta \, \beta' (\mathbf{v} - \mathbf{v}')\cdot  \Sigma^{-1}_g \cdot (\mathbf{v} - \mathbf{v}') f_g \, f_g' > 0
\end{equation}
for any $f$ with finite support, implying:
\begin{equation}
    n_\beta B_\beta - {\boldsymbol{A}}_\beta \cdot \boldsymbol{p}_\beta > 0.
\end{equation} 
The integral/infinite-sum exchange in \eqref{eq:den_pos_proof} is allowed for integrable positive integrands by the Tonelli-Fubini theorem.

\section{Derivation of the RFP Chang-Cooper weights}
\label{app:rfp-chang-cooper}
We start by deriving the RFP Chang-Cooper weight \eqref{modified_CC} in the perpendicular direction. The proof for the parallel direction follows straightforwardly. 
At equilibrium, the perpendicular flux vanishes, $J_{\perp,i+1/2,j}=0$, which yields
\begin{align}
    0
    &= D_{\perp\perp,i+1/2,j}\,\frac{f_{i+1,j}-f_{i,j}}{\Delta v_\perp}
       + a_{\perp,i+1/2,j}\Bigl[(1-\theta_{\perp,i+1/2,j})\,f_{i+1,j}
       + \theta_{\perp,i+1/2,j}\,f_{i,j}\Bigr] \\
    &= f_{i+1,j}\Bigl(
          \frac{D_{\perp\perp,i+1/2,j}}{\Delta v_\perp}
          + a_{\perp,i+1/2,j}(1-\theta_{\perp,i+1/2,j})
       \Bigr)
     + f_{i,j}\Bigl(
          -\frac{D_{\perp\perp,i+1/2,j}}{\Delta v_\perp}
          + a_{\perp,i+1/2,j}\,\theta_{\perp,i+1/2,j}
       \Bigr),
\end{align}
where $\theta_{\perp,i+1/2,j}$ is a yet-to-be-determined weight.  
Following the classical Chang-Cooper derivation, we determine $\theta_{\perp,i+1/2,j}$ by enforcing exact preservation of the steady state. From the previous relation,
\begin{equation}
    \frac{f_{i+1,j}}{f_{i,j}}
    =
    \frac{
        \displaystyle
        \frac{D_{\perp\perp,i+1/2,j}}{\Delta v_\perp}
        - a_{\perp,i+1/2,j}\,\theta_{\perp,i+1/2,j}
    }{
        \displaystyle
        \frac{D_{\perp\perp,i+1/2,j}}{\Delta v_\perp}
        + a_{\perp,i+1/2,j}\,(1-\theta_{\perp,i+1/2,j})
    }.
    \label{CC_ratio_numerical}
\end{equation}
In equilibrium, $f = f^M$ and we can write:
\begin{align}
    \frac{f^M_{i+1,j}}{f^M_{i,j}}
    &= \frac{
        \exp\!\Bigl(
            -\dfrac{v_{\perp,i+1}^2 + (v_{\parallel,j}-u_\parallel)^2}{v_{th,M}^2}
        \Bigr)
    }{
        \exp\!\Bigl(
            -\dfrac{v_{\perp,i}^2 + (v_{\parallel,j}-u_\parallel)^2}{v_{th,M}^2}
        \Bigr)
    } \\
    &= \exp\!\Bigl(
        -\dfrac{v_{\perp,i+1}^2 - v_{\perp,i}^2}{v_{th,M}^2}
    \Bigr)
     = \exp\!\Bigl(
        -\dfrac{(v_{\perp,i+1}+v_{\perp,i})\,\Delta v_\perp}{v_{th,M}^2}
    \Bigr) \\
    &= \exp\!\Bigl(
        -\,\Delta v_\perp\,\frac{v_{\perp,i+1/2}}{v_{th,M}^2}
    \Bigr)
     = \exp\!\bigl(-\Delta v_\perp\,w^M_{i+1/2}\bigr),    \label{CC_analytical}
\end{align}
where $v_{\perp,i+1/2} = \tfrac{1}{2}(v_{\perp,i+1}+v_{\perp,i})$ and
\[
    {w}^M_{i+1/2}
    := \frac{v_{\perp,i+1/2}}{v_{th,M}^2}.
\]
Equating the numerical ratio~\eqref{CC_ratio_numerical} with the analytical Maxwellian ratio~\eqref{CC_analytical}, and solving for the weight $\theta_\perp$, we obtain:
\begin{equation}
    \theta_{\perp,i+1/2,j}
    =
    \frac{1}{\Delta v_\perp w_{\perp,i+1/2,j}}
    - \frac{1}{\exp\!\bigl(\Delta v_\perp{w}^M_{\perp,i+1/2}\bigr)-1}.
\end{equation}
Similarly, by solving for the modified Chang-Cooper weight in the parallel direction, we obtain:
\begin{equation}
    \theta_{\parallel,i,j+1/2}
    =
    \frac{1}{\Delta v_\parallel w_{\parallel,i,j+1/2}}
    - \frac{1}{\exp\!\bigl(\Delta v_\parallel{w}^M_{\parallel,j+1/2}\bigr)-1}.
\end{equation}

% \appendix

\section{Discretization of Rosenbluth potentials Poisson equations}\label{Appendix:Poisson_Potentials}
At each time step, we compute the diffusion tensor $\mathcal{D}$ and the
advection vector $\mathbf{A}$ by solving the Poisson equations
\eqref{eq:rosenbluth_poisson}. 
We first obtain the potential $H$ from
\[
\nabla_\mathbf{v}^2 H = -\,8\pi f,
\]
where, in cylindrical velocity coordinates $(v_\perp, v_\parallel)$,
\[
\nabla_\mathbf{v}^2 H
= \frac{1}{v_\perp}\,\partial_{v_\perp}\!\bigl(v_\perp\,\partial_{v_\perp} H\bigr)
+ \partial_{v_\parallel}^2 H.
\] 

A second-order centered finite-difference discretization of
$\nabla_\mathbf{v}^2 H = -8\pi f$ then reads
\begin{align}
\frac{1}{v_{\perp,i}} \frac{1}{\Delta v_\perp}
\left[
  v_{\perp,i+\tfrac12}\,\frac{H_{i+1,j} - H_{i,j}}{\Delta v_\perp}
 -v_{\perp,i-\tfrac12}\,\frac{H_{i,j} - H_{i-1,j}}{\Delta v_\perp}
\right]
+\frac{H_{i,j+1} - 2H_{i,j} + H_{i,j-1}}{\Delta v_\parallel^2}
= -\,8\pi f_{i,j},
\label{eq:disc_poisson_H}
\end{align}
for $1 \le i \le N_\perp$ and $1 \le j \le N_\parallel$. The second Poisson equation,

\[
\nabla_\mathbf{v}^2 G = H,
\]
is discretized in an analogous manner. Denoting
$G_{i,j} \approx G(v_{\perp,i},v_{\parallel,j})$, we obtain
\begin{align}
\frac{1}{v_{\perp,i}} \frac{1}{\Delta v_\perp}
\left[
  v_{\perp,i+\tfrac12}\,\frac{G_{i+1,j} - G_{i,j}}{\Delta v_\perp}
 -v_{\perp,i-\tfrac12}\,\frac{G_{i,j} - G_{i-1,j}}{\Delta v_\perp}
\right]
+\frac{G_{i,j+1} - 2G_{i,j} + G_{i,j-1}}{\Delta v_\parallel^2}
= H_{i,j},
\label{eq:disc_poisson_G}
\end{align}
for $1 \le i \le N_\perp$ and $1 \le j \le N_\parallel$.

The ghost-cell values of $H$ and $G$ at the velocity-domain boundaries can
either be obtained from the Green's-function formalism
\cite{chacon2000implicit} or from a multipole-expansion approximation
\cite{taitano2016adaptive}. We adopt the latter for reduced computational
cost. In terms of the relative velocity
\[
\boldsymbol{c} = \boldsymbol{v} - \boldsymbol{u}, \qquad
c = |\boldsymbol{c}|,
\]
with cylindrical components $(c_\perp,c_\parallel)$, the {two-term-truncated} far-field
multipole expansions of the Rosenbluth potentials are 
\begin{align}
G(\boldsymbol{c})
&= n\,c 
 + \boldsymbol{\nabla}_{\boldsymbol{c}}\boldsymbol{\nabla}_{\boldsymbol{c}} c 
   : \left(\frac12 \int d\mathbf{c}'\, f(\boldsymbol{c}')\,\boldsymbol{c}'\boldsymbol{c}'\right) ,
\\[0.3em]
H(\boldsymbol{c})
&= \frac{2 n}{c}
 + \frac{3\,\boldsymbol{c}\boldsymbol{c}-\mathbf{I}\,c^{2}}{c^{5}}
   : \left(\int d\mathbf{c}'\, f(\boldsymbol{c}')\,\boldsymbol{c}'\boldsymbol{c}'\right) .
\end{align}

where $\mathbf{I}$ is the identity tensor and
$n = \int d\mathbf{c}'\, f(\boldsymbol{c}')$.

Using the identities \cite{taitano2016adaptive},
\begin{align}
\boldsymbol{\nabla}_{\boldsymbol{c}}\boldsymbol{\nabla}_{\boldsymbol{c}} c 
: \left(\frac12 \int d\mathbf{c}'\, f(\boldsymbol{c}')\,\boldsymbol{c}'\boldsymbol{c}'\right)
&= \frac{1}{c}\left(\int d\mathbf{c}'\, f(\boldsymbol{c}')\,\frac{c_\perp'^2}{2}\right)
 + \frac{c_\perp^{2}}{2 c^{3}}
   \left[\int d\mathbf{c}'\, f(\boldsymbol{c}')
          \left(c_\parallel'^2 - \frac12 c_\perp'^2\right)\right],
\\[0.3em]
\frac{3\,\boldsymbol{c}\boldsymbol{c}-\mathbf{I}\,c^{2}}{c^{5}}
: \left(\int d\mathbf{c}'\, f(\boldsymbol{c}')\,\boldsymbol{c}'\boldsymbol{c}'\right)
&= 2\,\frac{c_\parallel^{2}-\tfrac12 c_\perp^{2}}{c^{5}}
   \int d\mathbf{c}'\, f(\boldsymbol{c}')
        \left(c_\parallel'^2 - \frac12 c_\perp'^2\right).
\end{align}
We obtain the explicit far-field expressions
\begin{align}
G(\boldsymbol{c})
&= n\,c 
 + \frac{1}{c}\left(\int d\mathbf{c}'\, f(\boldsymbol{c}')\,\frac{c_\perp'^2}{2}\right)
 + \frac{c_\perp^{2}}{2 c^{3}}
   \left[\int d\mathbf{c}'\, f(\boldsymbol{c}')
          \left(c_\parallel'^2 - \frac12 c_\perp'^2\right)\right],
\\[0.3em]
H(\boldsymbol{c})
&= \frac{2 n}{c}
 + 2\,\frac{c_\parallel^{2}-\tfrac12 c_\perp^{2}}{c^{5}}
   \int d\mathbf{c}'\, f(\boldsymbol{c}')
        \left(c_\parallel'^2 - \frac12 c_\perp'^2\right),
\end{align}
which we use to populate the ghost cells for $H$ and $G$ at the outer
velocity-domain boundaries. In the discrete setting, the velocity integrals appearing in the
far-field formulas are approximated once per time step by quadrature
sums over the $(v_\perp,v_\parallel)$ grid. We use
$d\mathbf{c}' = 2\pi v_{\perp,i}\,\Delta v_\perp\,\Delta v_\parallel$ and set
\begin{align}
n
&:= \int d\mathbf{c}'\, f(\boldsymbol{c}')
 \;\approx\;
 2\pi \sum_{i=1}^{N_\perp}\sum_{j=1}^{N_\parallel}
 f_{i,j}\,v_{\perp,i}\,\Delta v_\perp\,\Delta v_\parallel,
\\[0.25em]
M_\perp
&:= \int d\mathbf{c}'\, f(\boldsymbol{c}')\,\frac{c_\perp'^2}{2}
 \;\approx\;
 2\pi \sum_{i=1}^{N_\perp}\sum_{j=1}^{N_\parallel}
 f_{i,j}\,\frac{c_{\perp,i}^{\,2}}{2}\,v_{\perp,i}\,\Delta v_\perp\,\Delta v_\parallel,
\\[0.25em]
Q
&:= \int d\mathbf{c}'\, f(\boldsymbol{c}')
      \Bigl(c_\parallel'^2 - \tfrac12 c_\perp'^2\Bigr)
 \;\approx\;
 2\pi \sum_{i=1}^{N_\perp}\sum_{j=1}^{N_\parallel}
 f_{i,j}\,\Bigl(c_{\parallel,j}^{\,2} - \tfrac12 c_{\perp,i}^{\,2}\Bigr)
 v_{\perp,i}\,\Delta v_\perp\,\Delta v_\parallel,
\end{align}
where $\boldsymbol{c}=(c_\perp,c_\parallel)$ is evaluated at the
grid points (e.g., $\boldsymbol{c}_{i,j}=\boldsymbol{v}_{i,j}-\boldsymbol{u}$).
Using the discrete moments $n$, $M_\perp$, and $Q$, the ghost-cell
values of the Rosenbluth potentials are then given by the discrete
far-field formulas
\begin{align}
 H_{i,j}^{\text{ghost}}
&= \frac{2 n}{c_{i,j}}
 + 2\,\frac{c_{\parallel,j}^{2}-\tfrac12 c_{\perp,i}^{2}}{c_{i,j}^{5}}\,Q,
\\[0.25em]
G_{i,j}^{\text{ghost}}
&= n\,c_{i,j}
 + \frac{M_\perp}{c_{i,j}}
 + \frac{c_{\perp,i}^{2}}{2\,c_{i,j}^{3}}\,Q,
\end{align}
for the boundary indices
\[
(i,j)\in\{0,N_\perp+1\}\times\{1,\dots,N_\parallel\}
\quad\text{or}\quad
(i,j)\in\{1,\dots,N_\perp\}\times\{0,N_\parallel+1\}.
\]
These ghost cell values are then used in the discrete Poisson solves
\eqref{eq:disc_poisson_H}-\eqref{eq:disc_poisson_G}.
Following this, the diagonal diffusion tensor components are obtained by a
second-order centered finite-difference approximation of $G$,
\begin{align}
D_{\perp\perp}\big|_{i,j}
&:= \partial_{c_{\perp}}^{2} G(\boldsymbol{c}_{i,j})
 \;\approx\;
 \frac{G_{i+1,j}-2G_{i,j}+G_{i-1,j}}{\Delta v_\perp^{2}}, 
\\[0.3em]
D_{\parallel\parallel}\big|_{i,j}
&:= \partial_{c_{\parallel}}^{2} G(\boldsymbol{c}_{i,j})
 \;\approx\;
 \frac{G_{i,j+1}-2G_{i,j}+G_{i,j-1}}{\Delta v_\parallel^{2}}.
\end{align}
For the off-diagonal tensor component, we use the standard second-order,
centered expression
\begin{equation}
D_{\perp\parallel}\big|_{i,j}
\;\approx\;
\frac{
  G_{i+1,j+1} - G_{i+1,j-1}
 -G_{i-1,j+1} + G_{i-1,j-1}}
{4\,\Delta v_\perp\,\Delta v_\parallel}.
\end{equation}
Finally, the advection vector $\mathbf{A} = (A_\perp,A_\parallel)$ is
computed directly at the cell faces from $H_{i,j}$ as:
\begin{align}
A_{\perp,i+\tfrac12,j}
&= \partial_{v_\perp} H\big|_{i+\tfrac12,j}
 \;\approx\;
 \frac{H_{i+1,j} - H_{i,j}}{\Delta v_\perp},
\\[0.3em]
A_{\parallel,i,j+\tfrac12}
&= \partial_{v_\parallel} H\big|_{i,j+\tfrac12}
 \;\approx\;
 \frac{H_{i,j+1} - H_{i,j}}{\Delta v_\parallel}.
\end{align}

\bibliographystyle{elsarticle-num}
\bibliography{Paper_Draft/references}

\end{document}